\documentclass[12pt]{amsart}
\usepackage{amsmath,amsfonts,amsthm,amssymb}
\newcommand\R{{\mathbb R}}

 \newtheorem{theorem}{Theorem}
 \newtheorem{conjecture}{Conjecture}
 \newtheorem{remark}{Remark}[section]
 \newtheorem{lemma}[remark]{Lemma}
 \newtheorem{cor}[remark]{Corollary}
 \newtheorem{proposition}[remark]{Proposition}
 \newtheorem{definition}[remark]{Definition}
\newcommand\e{{\epsilon}}

\begin{document}
\title{A-priori bounds for the 1-d cubic NLS in negative Sobolev spaces}

\author{Herbert Koch}
\address{Mathematisches Institut   \\ Universit\"at Bonn }

\author{ Daniel Tataru}
\address {Department of Mathematics \\
  University of California, Berkeley} \thanks{ The first author was
  partially supported by DFG grant KO1307/1, by MSRI for Fall
  2005 and by the Miller Institute for basic research in Science in Spring 2006\\ The second author was partially supported by NSF grants
  DMS0354539 and DMS 0301122 and also by MSRI for Fall 2005}

\begin{abstract}
  We consider the cubic Nonlinear Schr\"odinger Equation (NLS)  in one space dimension, either
  focusing or defocusing.  We prove that the solutions satisfy
  a-priori local in time $H^{s}$ bounds in terms of the $H^s$ size of
  the initial data for $s \geq -\frac16$.
\end{abstract}
\maketitle

\section{Introduction}

The one dimensional cubic Nonlinear Schr\"odinger  equation (NLS) 
\begin{equation}
i u_t - u_{xx} \pm u |u|^2 = 0, \qquad u(0) = u_0.
\label{nls}\end{equation}
arises as generic asymptotic equation for modulated wave trains.
  Its has a particularly rich structure: It is Hamiltonian with respect to 
the symplectic structure  
\[   \sigma(u,v) = \text{Im} \int u\overline{v} dx    \]
and the Hamiltonian
\[  \int \frac12(u')^2 \pm \frac14 |u|^4 dx. \]
There are infinitely many conserved quantities.  The NLS equation is
completely integrable in the sense that there exist Lax pairs for
it. The machinery of inverse scattering allows to construct many
interesting solutions, among them solitary waves in the focusing
case.
 
The NLS  is globally well-posed for initial data $u_0 \in L^2$, and 
locally in time the solution has a uniform lipschitz  dependence on the
initial data in balls.

On the other hand \eqref{nls} is invariant with respect to the scaling
\[
u(x,t) \to \lambda u(\lambda x,\lambda^2 t) 
\]
This implies that the scale invariant initial data space for
\eqref{nls} is $\dot H^{-\frac12}$. Thus one is motivated to ask
whether the local well-posedness also holds in negative Sobolev
spaces. 

The equation \eqref{nls} is also invariant under the Galilean
transformation
\[ 
u(x,t) \to   e^{icx-ic^2 t}     u (x+2ct,t) 
\]
which corresponds to a shift in the frequency space.  As a consequence
there is no uniformly continuous dependence on the initial data (see \cite{MR1813239},
\cite{MR2018661}). This is not unexpected; if local uniformly continuous
dependence were to hold in any negative Sobolev space, by Gallilean
invariance and scaling this would imply global in time local in space
uniformly continuous dependence  on the initial data in $L^2$.

What we expect below $L^2$ is for the cubic NLS to exhibit genuinely
nonlinear dynamics, which corresponds to a continuous but not
uniformly continuous dependence on the initial data. One may be
tempted to think that local well-posedness should hold all the way
down to $s = -\frac12$.  However, such a result is far out of reach
for now and we would not even speculate whether it is true or not.

On the other hand, there is another very natural threshold, which 
is connected to the main motivation of the present paper. In a recent
paper Kappeler and Topalov~\cite{MR2131061}  proved that the mKdV equation 
\[
v_t - v_{xxx}+ v_x v^2 = 0, \qquad v_0 = v_0
\]
on the torus is well-posed for initial data in $L^2$. The proof relies
on complete integrability of the equation, and it uses the machinary
of integrable equations in a fundamental way.
One may ask whether the
same result holds on the real line, and also whether it is possible to
find  arguments  which do not use the integrable structure.

To connect this problem with the NLS equation we consider modulated
wave train solutions $v$ of the form $v = \Re w$ where $w$ is
frequency localized in a neighborhood of size $h$ of some large
frequency $\lambda$.  Then $w$ solves the equation
\[
w_t +i \lambda^3 w + 3i \lambda^2(D_x-\lambda) w + 3i
\lambda(D_x-\lambda)^2 w + 3 i \lambda w |w|^2 \! \approx O(h^3) w +
O(h) w|w|^2 + w^3
\]
For $h \ll \lambda$ we neglect the first two terms on the right. The
$w^3$ term is non-resonant and is also neglected. Then the substitution
\[
w(t,x)  = \lambda^{-\frac12} e^{-i \lambda^3 t} e^{i\lambda x}
u(t,\lambda^{-\frac12}(x-3\lambda^2 t)) 
\]
turns the above equation into \eqref{nls} with modified constants.

A frequency range of size $\lambda$ for $u$ turns into a frequency
range of size $\mu = \lambda^\frac32$ for $w$.  By construction this
frequency range for $w$ is centered at the origin, but we can use a
Galilean transformation to shift it to a dyadic region.  We can also
easily compute
\[
\| u(0)\|_{L^2} = \lambda^{\frac14} \|v(0)\|_{L^2} = 
\mu^\frac16  \|v(0)\|_{L^2}
\]
Hence the mKdV equation with initial data in $L^2$ is
similar\footnote{We emphasize that this similarity applies only 
for solutions in a dyadic frequency range. On the other hand
in our analysis later in the paper we see that some of the 
most difficult to control multilinear interactions occur in the case
of unbalanced frequencies, where this analogy no longer applies.}
 to the NLS
equation with initial data in $H^{-\frac16}$.  We view the one
dimensional NLS equation as a simpler model in the analysis of the KdV
equation; this is due to the added Gallilean invariance. However, it
is also interesting in its own right.

The threshold $s=-\frac16$ also arises in several key steps of our
analysis later on, having to do with the interaction of high and low 
frequencies.  We are led to

\begin{conjecture}
The cubic NLS equation \eqref{nls} is locally well-posed for initial
data in $H^s$ with $s \geq -\frac16$.
\end{conjecture}

To prove this one would need to establish a-priori $H^s$ bounds for the
solutions and then prove continuous dependence on the initial data.
In this article we solve the easier half of this problem.

\begin{theorem}\label{main} 
  Let $s \geq -\frac16$. For any $M > 0$ there exists $T > 0$ and $C >
  0$ so that for any initial data $u_0 \in L^2$ satisfying
\[
\| u_0\|_{H^s} \leq M
\]
there exists a solution $u \in C(0,T;L^2)$ to \eqref{nls} 
which satisfies 
\[
\| u\|_{L^\infty H^{s}} \leq C \| u_0\|_{H^s}
\]
\end{theorem}

While writing  this paper the authors have learned that similar
results were independently obtained by Christ-Colliander-Tao~\cite{CCT}. 
Their results apply in the range $s > -\frac1{12}$. 

We also refer the reader to the work of Vargas-Vega~\cite{MR1876762}
and Gr{\"u}nrock~\cite{MR2181058} who consider the cubic NLS in
alternative function spaces below $L^2$, but only in settings where the
local Lipschitz dependence on the initial data still holds.
 
\begin{remark}
  In the process of proving the theorem we actually obtain a better
  characterization of the solution $u$, namely we show that $u$
  bounded in a space $X^s$ defined in the next section which embeds
  into $L^\infty H^s$ and has the property that the nonlinear
  expression $|u|^2u$ is well defined for $u \in X^s$ with a bound 
 depending only on the $H^s$ norm of the initial data. 
\end{remark}

We note that by rescaling the problem reduces to the case of small
initial data. Then we take $M=\e$, small and $T=1$, $C=2$.

We begin with a dyadic frequency decomposition of  the solution $u$,
\[
u = \sum_\lambda u_\lambda
\]
To measure the $H^s$ norm of $u$ we use  the stronger norm than $L^\infty(H^s)$,
\[
\|u\|_{l^2 L^\infty H^s}^2 =  \sum_\lambda \sup_t \lambda^{2s} \| u_\lambda
(t)\|_{L^2}^2 
\]
That we can use this instead of the $L^\infty H^s$ norm is a
reflection of the fact that there is not much energy transfer between
different dyadic frequencies.

To prove the theorem we need two Banach spaces $X^s$ and $Y^s$,
defined in the next section, in order to measure the regularity of the
solution $u$, respectively of the nonlinear term $|u|^2 u$.

The linear part of the argument is given by 
\begin{proposition}
The following estimate holds:
\[
\| u\|_{X^s} \lesssim  \| u \|_{l^2L^\infty H^{s}} + \| iu_t -\Delta u\|_{Y^s}
\]
\label{plin}\end{proposition}
To estimate the nonlinearity we need a cubic bound,

\begin{proposition}
Let $-\frac16 \leq s \leq 0$ and $u \in X^s$. Then $|u|^2 u \in Y^s$ and
\[
\| |u|^2 u\|_{Y^s} \lesssim \|u\|_{X^s}^3
\]
\label{pnonlin}\end{proposition}

Finally we need to propagate the $H^s$ norm:

\begin{proposition}
Let $-\frac16 \leq s \leq 0$, and $u$ be a solution to \eqref{nls} with
\[ 
\| u\|_{l^2 L^\infty H^{s}} \ll 1. 
\]
Then  we have
\[
\| u\|_{l^2 L^\infty H^{s}} \lesssim \| u_0\|_{H^{s}} 
  + \|u\|_{X^s}^3.
\]
\label{penergy}\end{proposition}

The plan of the paper is as follows. In the next section we motivate 
and introduce the spaces $X^s$ and $Y^s$, as well as establish the
linear mapping properties in Proposition~\ref{plin}. In Section~\ref{sbilin}
we discuss the linear and bilinear Strichartz estimates for solutions
to the linear equation.  

The trilinear estimate in Proposition~\ref{pnonlin} is proved in
Section~\ref{snonlin}. Finally in the last section we use a variation
of the I-method to construct a quasi-conserved energy functional and
compute its behavior along the flow, thus proving
Proposition~\ref{penergy}.

To conclude this section we show that the conclusion of the Theorem
follows from the above Propositions. We first note that if  $u_0 \in
L^2$ then by iteratively solving the equation on small time intervals
we obtain a solution $u$ up to time $1$, which satisfies
\begin{equation}
i u_t -\Delta u \in L^2
\label{pk}\end{equation}
This easily implies that $u \in l^2 L^\infty H^{s}$, and also that $u
\in X^s$. 

To prove the theorem we use a continuity argument. Let $\varepsilon>0$ be a
small constant and suppose that $\Vert u_0 \Vert_{H^s(\R)} < \varepsilon$.
Fix a small threshold $\delta$, $\varepsilon \ll \delta \ll 1 $ and denote by
$A$ the set
\[
A = \{ T \in [0,1]; \  \|u\|_{l^2 L^\infty H^{s}([0,T] \times \R)} \leq 2
\delta,\ \|u\|_{X^s([0,T] \times \R)} \leq 2 \delta \}
\]
We claim that $A = [0,1]$. To show this we first observe that $0 \in
A$. The norms above increase with $T$, therefore $A$ is an interval.
We show that $A$ is both open and closed in $[0,1]$. 

By \eqref{pk} it easily follows that the norms in the definition of
$A$ are continuous~\footnote{This of course depends on the definition
  of the $X^s$ norm, but it is straightforward to prove.} with respect
to $T$. This implies that $A$ is closed.

Finally let $T \in A$. By Proposition~\ref{penergy} we obtain 
\[
\|u\|_{l^2 L^\infty H^{s}([0,T] \times \R)} \lesssim \varepsilon + 
\delta^3.
\]
Then by Propositions~\ref{plin}, \ref{pnonlin} we obtain 
\[
\|u\|_{X^s([0,T] \times \R)} \lesssim \e + \delta^3 
\]
If $\e$ and $\delta$ are chosen to be sufficiently small we conclude
that 
\[
\|u\|_{l^2 L^\infty H^{s}([0,T] \times \R)} \leq \delta, \qquad
\|u\|_{X^s([0,T] \times \R)} \leq \delta
\]
By the continuity of the norms with respect to $T$ it follows that a
neighborhood of $T$ is in $A$. 

Hence $A = [0,1]$ and the Theorem \ref{main} is proved.

\section{ The function spaces}
To understand what to expect in terms of the regularity of $u$ we
begin with some heuristic considerations.
If the initial data $u_0$ to \eqref{nls} satisfies $\|u_0\|_{L^2} \leq
1$ then the equation can be solved iteratively using the Strichartz
estimates. We obtain essentially linear dynamics, and the solution $u$
belongs to the space $X^{0,1}$ associated to the Schr\"odinger
equation (see the definition in \eqref{Xsb} below). 

Let $s < 0$. Consider now the same problem but with initial data $u_0 \in H^{s}$,
localized at frequency $\lambda$. Then the initial data satisfies
$\|u_0\|_{L^2} \lesssim \lambda^{-s}$. By rescaling we conclude that
the evolution is still described by linear dynamics up to the time
$\lambda^{4s}$. 

Then it is natural to consider a dyadic decomposition of  the solution $u$ 
\[
u = \sum_\lambda u_\lambda
\]
and to measure the $u_\lambda$ component uniformly in $\lambda^{4s}$
time intervals. We remark that this is reasonable for as long as there is
not much input coming from the higher frequencies. This is the technical point 
 where the $s= -\frac16$ threshold arises in our proof. 

A good candidate for measuring  $u_\lambda$ in $\lambda^{4s}$
time intervals is given by Bourgain's  $X^{s,b}$ spaces defined by
\begin{equation} \label{Xsb}
\| u\|_{X^{s,b}}^2 = \int |\hat{u}(\tau,\xi)|^2 \xi^{2s} (1+
|\tau-\xi^2|)^{2b} d\xi d\tau
\end{equation} 
where the natural choice for $b$ from a scaling standpoint is
$b=\frac12$. However, this choice leads to logarithmic divergences
in estimates, so one commonly uses instead some $b > \frac12$ but
close to it. We could do this here but it would complicate the
bookkeeping and would also not work at $s = -\frac16$.
For $b=\frac12$ one can go one step further and consider dyadic
decompositions with respect to the  modulation $\tau-\xi^2$.
This leads to  the additional homogeneous Besov type norms
\[
\| u\|_{\dot{X}^{s,\frac12,1}} = \sum_\mu \left( \int_{|\tau-\xi^2| \approx
    \mu} |\hat{u}(\tau,\xi)|^2 \xi^{2s} |\tau-\xi^2| d\xi
  d\tau\right)^\frac12
\]
\[
\| u\|_{\dot{X}^{s,\frac12,\infty}} = \sup_\mu \left( \int_{|\tau-\xi^2| \approx
    \mu} |\hat{u}(\tau,\xi)|^2 \xi^{2s} |\tau-\xi^2| d\xi
  d\tau\right)^\frac12
\]

Instead in this paper we use the closely related spaces $U^2_\Delta$
and $V^2_\Delta$.  Spaces of this type have been first introduced in
unpublished work of the second author on wave-maps, but in the
meantime they have been also used in \cite{MR2094851}, \cite{BeTa}, \cite{SS}. They
turn out to be useful replacements of $X^{s,b}$ spaces in limiting
cases, and they retain the scaling of the corresponding space of
homogeneous solutions to the linear equation.  We define them and
summarize their key properties in what follows.

\begin{definition}
  Let $1 \leq p < \infty$. Then $U^p_\Delta$ is an atomic space, where
  atoms are piecewise solutions to the linear equation,
\[
u = \sum_k 1_{[t_k,t_{k+1})} e^{it D_x^2} u_k, \qquad \sum_k
\|u_k\|_{L^2}^p = 1
\]
and $\{t_k\}$ is an arbitrary increasing sequence. 
\end{definition}
Clearly we have
\[
U^p_\Delta \subset L^\infty L^2
\]
In addition, the $U^p_\Delta$ functions are continuous except at
countably many points, and right continuous everywhere.

A close relative is the space $V^p_\Delta$ of functions with bounded
$p$-variation along the flow:

\begin{definition}
  Let $1 \leq p < \infty$. Then $V^p_\Delta$ is the space of right
  continuous functions $u \in L^\infty(L^2)$ for which the following
  norm is finite,
\[
\|u\|_{V^p_\Delta}^p = \|u\|_{L^\infty L^2}^p+\sup_{\{t_k\}\nearrow}
\sum_k \| e^{i t_k D_x^2} u(t_{k}) - e^{i t_{k+1} D_x^2 }u(t_{k+1})\|_{L^2}
\] 
where the supremum is taken with respect to all increasing sequences
$\{t_k\}$.
\end{definition}

Conjugation with the Schr\"odinger group reduces a large part of the study 
of the spaces $V^p$ and $U^p$ to the scalar case, where we replace the 
group by the identity.

We have the series of inclusions
\begin{equation}
 U^p_\Delta \subset V^p_\Delta \subset
U^q_\Delta \subset
 L^\infty L^2, \qquad p < q.
\label{embedding}\end{equation} 
The inclusion $U^p \subset V^p$ can easily checked on atoms. The imbedding 
$V^p \subset U^q$ is a little harder and its proof can be found in 
Section 5 of \cite{MR2094851}.

We denote by $DU^p_\Delta$ the space of functions
\[
DU^p_\Delta = \{ (i\partial_t -\partial_x^2) u; \ u \in U^p_\Delta \}
\]
with the induced norm.  Then we have the trivial bound
\begin{equation}
\| u\|_{U^p_\Delta} \lesssim 
\| u(0)\|_{L^2} + \| (i\partial_t -\partial_x^2) u\|_{DU^p_\Delta}
\label{lin}\end{equation}
Finally, we have the duality relation
\begin{equation}
(DU^p_\Delta)^* = V^{p'}_\Delta. 
\label{dual}\end{equation}
To see this one first verifies the inequality 
\[ | \int \langle(i\partial_t - \partial_x^2)f,g\rangle_{L^2}  dx | 
\le  \Vert f \Vert_{U^p} \Vert g \Vert_{V^{p'}}   \]
by checking it  for atoms $f$. Secondly, given  $L \in (DU^p_\Delta)^*$ 
we apply it to characteristic functions of  intervals, which allows to define a function 
$g$ with 
\[ \int \langle(i\partial_t - \partial_x^2) f , g \rangle_{L^2} dt =
L( (i\partial_t - \partial_x^2)f ). \] An application to suitable
atoms shows that $g \in V^{p'}_\Delta$.

Moreover we have the embedding
\[ \dot X^{0, \frac12,1} \subset U^2_\Delta.  \]
To see this it suffices to consider a function  the Fourier transform of which
is supported in a fixed dyadic annulus. The statement follows now
easily. Combined with duality one sees that 

\begin{equation} \label{embeddingX} 
\dot X^{0, \frac12,1} \subset U^2_\Delta \subset V^2_\Delta \subset \dot
X^{0,\frac12,\infty}. 
\end{equation}

The $U^p_\Delta$ and $V^p_\Delta$ spaces behave well with respect to
sharp time truncations.  Precisely, if $I$ is a time interval and
$\chi_I$ is its characteristic function then we have the
multiplicative mapping properties
\begin{equation}
  \chi_I:  U^p_\Delta \to U^p_\Delta, \qquad 
\chi_I:  V^p_\Delta \to V^p_\Delta
\label{mult}\end{equation}
with uniform bounds with respect to $I$.

We use a spatial Littlewood-Paley decomposition
\[
1 = \sum_{\lambda \geq 1\ dyadic} P_\lambda, \qquad u 
= \sum_{\lambda  \geq 1\ dyadic} P_\lambda u 
= \sum_{\lambda \geq 1\ dyadic}
u_\lambda
\]
as well as a  Littlewood-Paley decomposition with respect to the 
modulation $\tau-\xi^2$,
\[
1 = \sum_{\lambda \geq 1\ dyadic}  Q_\lambda
\]
Both decompositions are inhomogeneous. It is easy to verify 
that we have the uniform boundedness properties
\begin{equation}
P_\lambda:  U^p_\Delta \to U^p_\Delta, \qquad Q_\lambda: U^p_\Delta
\to U^p_\Delta
\label{pqbd}\end{equation}
and similarly for $V^p_\Delta$.

For functions at frequency $\lambda$ we introduce a minor variation
of the $U^2_\Delta$, respectively $V^2_\Delta$ spaces, which we 
denote by $U^2_\lambda$, respectively $V^2_\lambda$. Their norms are
defined as
\[
\| u_\lambda \|_{ U^2_\lambda}^2 = \| Q_{\leq \lambda^2} u_\lambda\|_{
  U^2_\Delta}^2 + \sum_{|I| = \lambda^{-2}, |J| = \lambda^{-1}} 
\| \chi_I(t) \chi_J(x) Q_{\geq \lambda^2} u_\lambda \|_{U^2}^2,
\]
respectively
\[
\| u_\lambda \|_{ V^2_\lambda}^2 = \| Q_{\leq \lambda^2} u_\lambda\|_{
  V^2_\Delta}^2 + \sum_{|I| = \lambda^{-2}, |J| = \lambda^{-1}} 
\| \chi_I(t) \chi_J(x) Q_{\geq \lambda^2} u_\lambda \|_{V^2}^2,
\]
Here the time truncation is still sharp, as above. The spatial
truncation may be taken sharp or smooth, the two norms are equivalent
due to the frequency localization. In the last norm we use the simpler
space $U^2$ (where we replace $\Delta$ in $U_\Delta$ by zero) instead
of $U^2_\Delta$; this is also immaterial, the $U^2$ and $U^2_\Delta$
norms are equivalent at frequency $\lambda$ and modulation $\geq
\lambda^2$.

In doing this the $U^2_\Delta$ norm is slightly weakened, but only in the elliptic region:
\[
\| u_\lambda \|_{U^2_\lambda} \lesssim \|u_\lambda\|_{U^2_\Delta}
\]
To see this it suffices to consider $U^2_\Delta$ atoms. Steps
$t_{k+1}-t_k$ of size larger than $\lambda^{-2}$ are essentially
canceled by the modulation localization operator $Q_{\geq \lambda^2}$,
therefore is suffices to restrict ourselves to the $\lambda^{-2}$ time
scale. But on this scale the Schr\"odinger flow at frequency $\lambda$
is trivial, i.e. there is no propagation. Thus one obtains the square
summability with respect to the $\lambda^{-1}$ spatial scale.

Since we preserve the duality relation \eqref{dual} the $V^2_\Delta$
norm is slightly strengthened:
\[
\|  u_\lambda \|_{V^2_\lambda} \gtrsim \|u_\lambda\|_{V^2_\Delta}
\]

The only advantage in using the modified spaces is that they allow us 
to replace a logarithm of the high frequency by a logarithm of the low
frequency in \eqref{biuv}, which is needed in order for our 
proofs to work in the limiting case $s = -\frac16$.

We note that the inclusions in \eqref{embedding} as well as the
properties \eqref{lin}, \eqref{dual}, \eqref{embeddingX}, \eqref{mult}
and \eqref{pqbd} remain valid in the dyadic setting for the modified 
spaces.

Now we are ready to introduce the function spaces for the solutions
$u$. We set
\begin{equation}
\| u\|_{X^s}^2 = \sum_{\lambda} \lambda^{2s}  \sup_{|I| = \lambda^{4s}}
\|\chi_I u_\lambda\|_{U^2_\lambda}^2  
\label{xs}\end{equation}
where we sum over all dyadic integers $\ge 1$ with the obvious modification 
at $\lambda =1$.

To measure the regularity of the nonlinear term we need
\begin{equation}
\| f\|_{Y^s}^2 = \sum_{\lambda} \lambda^{2s}   \sup_{|I| = \lambda^{4s}}
\|\chi_I f_\lambda\|_{DU^2_\lambda}^2  
\label{ys}\end{equation}
Due to \eqref{lin} we easily obtain the bound in Proposition~\ref{plin}.

\section{Linear and bilinear estimates}
\label{sbilin}

We begin with solutions to the homogeneous equation,
\begin{equation}
iv_t -\Delta v = 0, \qquad v(0) = v_0
\label{hom}\end{equation}
These satisfy the Strichartz estimates:

\begin{proposition}
Let $p,q$ be indices satisfying 
\begin{equation}
\frac{2}p + \frac{1}q = \frac12, \qquad 4 \leq p \leq \infty
\label{pq}\end{equation}
Then the solution $u$ to \eqref{hom} satisfies 
\[
\|v\|_{L^p_t L^q_x} \lesssim \|v_0\|_{L^2}
\]
\end{proposition}

In particular we note the pairs of indices $(\infty,2)$, $(6,6)$ and
$(4,\infty)$.  On occasion it is convenient to interchange the role of
the space and time coordinates. Then by interpolating the 
local smoothing estimate for solutions to \eqref{hom}, 
\[
 \Vert v_\lambda \Vert_{L^\infty_xL^2_t} \lesssim \lambda^{-1/2}
 \Vert v_{\lambda,0} \Vert_{L^2} 
\]
and the maximal function estimate 
\[
 \Vert v_\lambda \Vert_{L^4_xL^\infty_t} \lesssim \lambda^{1/4}
 \Vert v_{\lambda,0} \Vert_{L^2} 
\] 
we obtain

\begin{proposition}
Let $p,q$ be indices satisfying 
\eqref{pq}.
Then for every solution $v$ to \eqref{hom} which is localized 
at frequency $\lambda$ we have 
\[
\|v_\lambda\|_{L^p_x L^q_t} \lesssim \lambda^{\frac{3}p-\frac12}  \|v_{\lambda,0}\|_{L^2}
\]
\end{proposition}

As a straightforward consequence we have

\begin{cor}
a) Let $p,q$ be indices satisfying  \eqref{pq}. 
Then 
\[
\|v\|_{L^p_t L^q_x} \lesssim \|v\|_{U^p_\Delta}
\]
and the same holds  with $U^p_\Delta $ replaced by $V^2_\Delta$.

b) In addition, if  $v$ is  is localized 
at frequency $\lambda$ then we have 
\[
\|v\|_{L^p_x L^q_t} \lesssim
\lambda^{\frac{3}p-\frac12}\|v\|_{U^p_\Delta}, 
\]
and the same holds  with $U^p_\Delta $ replaced by $V^2_\Delta$ if $p
> 2$.

c) For $v$ localized at frequency $\lambda$ the $U^p_\Delta $ 
and  $V^2_\Delta$ norms in (a), (b) can be replaced by $U^2_\lambda$ 
and  $V^2_\lambda$.
\label{upstr}\end{cor}
The proof is straightforward, since it suffices to do it for atoms. In
the case of $V^2_\Delta$ we also take advantage of the inclusion 
$V^2_\Delta \subset U^p_\Delta$, $p > 2$. The estimate for $V^2_\lambda$ and
$U^2_\lambda$  follows from the embeddings 
$U^2_\lambda \subset V^2_\lambda \subset V^2_{\Delta}$ for functions at
frequency $\lambda$.

By duality we also obtain

\begin{cor}
a) Let $p,q$ be indices satisfying  \eqref{pq}.
Then 
\[
 \|v\|_{DU^{2}_\Delta}  \lesssim \|v\|_{L^{p'}_t L^{q'}_x} 
\]

b) In addition, if  $v$ is  is localized 
at frequency $\lambda$  then we have 
\[
\|v\|_{DU^2_\Delta}  \lesssim
\lambda^{\frac{3}p-\frac12} \|v\|_{L^{p'}_x L^{q'}_t}, \qquad \text{ for } p > 2.
\]
c) For $v$ localized at frequency $\lambda$ the $DU^2_\Delta $ 
 norm in (a), (b) can be replaced by $DU^2_\lambda$.
\label{lpdual}\end{cor}

The second type of estimates we use are bilinear:

\begin{proposition}\label{bilinear} 
Let $\lambda > 0$. Assume that $u,v$ are solutions to \eqref{hom}
which are $\lambda$ separated in frequency. Then 
\begin{equation}
\| u v\|_{L^2} \lesssim \lambda^{-\frac12} \|u_0\|_{L^2} \|v_0\|_{L^2}
\end{equation}
\end{proposition}

\begin{proof}
In the Fourier space we have 
\[
\hat u(\tau,\xi) = \hat u_0(\xi) \delta_{\tau - \xi^2}, \qquad
\hat v(\tau,\xi) = \hat v_0(\xi) \delta_{\tau - \xi^2}
\]
Then 
\[
\widehat{(uv)}(\tau,\xi) = \int_{\xi_1+\xi_2 = \xi}  
 \hat u_0(\xi_1) \hat v_0(\xi_2) \delta_{\tau - \xi_1^2-\xi_2^2} d \xi_1
\]
which gives
\[
\widehat{(uv)}(\tau,\xi) =  \frac1{2|\xi_1-\xi_2|} (\hat u_0(\xi_1) \hat v_0(\xi_2) 
+ \hat  u_0(\xi_2) \hat v_0(\xi_1))
\]
where $\xi_1$ and $\xi_2$ are the solutions to
\[
\xi_1^2 +\xi_2^2 = \tau, \qquad \xi_1+\xi_2 = \xi
\]
We have 
\[
d\tau d\xi = 2|\xi_1-\xi_2| d\xi_1 d\xi_2
\]
therefore we obtain
\[
\| uv\|_{L^2} \lesssim  \int |\hat u_0(\xi_1)|^2 |\hat v_0(\xi_2)|^2
|\xi_1-\xi_2|^{-1}d\xi_1 d\xi_2
\]
The conclusion follows.
\end{proof}

As a consequence we obtain

\begin{cor}
a) Let $u,v$ be functions which are $\lambda$ separated in
frequency. Then
\begin{equation}
  \| u v\|_{L^2} \lesssim \lambda^{-\frac12} \|u\|_{U^2_\Delta} \|v\|_{U^2_\Delta}
\end{equation}

b) Let $\lambda \ll \mu$. Then 
\begin{equation}
  \| u_\lambda  v_\mu\|_{L^2} \lesssim \mu^{-\frac12} 
  \|u_\lambda \|_{U^2_\lambda} \|v_\mu\|_{U^2_\mu}
\label{biuu}\end{equation}
\end{cor}
Again it suffices to prove these estimates for atoms, and then for solutions to
the homogeneous Schr\"oder equation. But this follows from the $L^4$ Strichartz
estimates and the bilinear estimate of Proposition \ref{bilinear}.

At a single point in the paper we need a version of \eqref{biuu} with 
$U^2_\lambda$ replaced by $V^2_\lambda$. This is the only place where
we use the $V^2_\lambda$ modification of $V^2_\Delta$.

\begin{proposition}
Let $\lambda \ll \mu$ and $|I| = 1$. Then 
\begin{equation}
  \| \chi_I u_\lambda  v_\mu\|_{L^2} \lesssim \mu^{-\frac12} \ln \lambda 
\|u_\lambda \|_{V^2_\lambda} \|v_\mu\|_{U^2_\mu}
\label{biuv}\end{equation}
\end{proposition}

We note that in order to treat the limiting case $s=-\frac16$ it is
acceptable to loose $\ln \lambda$, but not $\ln \mu$.

\begin{proof}
We split $u_\lambda$ into a low modulation part and a high modulation
part, 
\[
u_\lambda = Q_{\leq \lambda^2} u_\lambda + Q_{\geq \lambda^2} u_\lambda 
\]
The first term is estimated in $V^2_\Delta$ simply by counting dyadic
regions with respect to modulation. The time truncation regularizes
the modulation less than $1$, so we are left with about $\log \lambda$
regions. 

On the other hand for the second term we use the $l^2$ summability
with respect to rectangles of size $\lambda^{-2} \times \lambda^{-1}$.
Precisely, via Bernstein's inequality we have
\[
\| Q_{\geq \lambda^2} u_\lambda\|_{l^2 L^\infty} \lesssim \lambda^{\frac12} \|u_\lambda\|_{V_\lambda^2}
\]
It remains to show that
\[
\|v_\mu\|_{L^2(R)} \lesssim \lambda^{-\frac12} \mu^{-\frac12} \|v_\mu\|_{U^2_\mu}
\]
where $R$ is a rectangle as above.  By the definition of $U^2_\mu$ the
problem reduces to the case when $v_\mu$ solves the homogeneous
Schr\"odinger equation. But in that case the above inequality is
nothing but the classical local smoothing estimate.

\end{proof}

\section{The cubic nonlinearity}
\label{snonlin}
In this section we prove Proposition~\ref{pnonlin}.  

For a dyadic frequency $\lambda$ we estimate the nonlinearity $|u|^2
u$ at frequency $\lambda$ in a $\lambda^{4s}$ time interval $I$.  We
take a dyadic decomposition of each of the factors and denote the
corresponding frequencies by $\lambda_1$, $\lambda_2$, $\lambda_3$. We
consider several cases:

{\bf Case 1.} $\lambda_{1,2,3} \lesssim \lambda$.
Then the $X^s$ bounds at the $\lambda_j$ frequencies are 
localized to time intervals at least as large as  $I$.
 Hence we use directly the $L^6$ Strichartz estimates to obtain
\[
\begin{split}
\lambda^s \| \chi_I u_{\lambda_1}  \bar{u}_{\lambda_2} u_{\lambda_3} 
\|_{DU^2_\lambda} & \lesssim \lambda^{3s} 
 \| \chi_I u_{\lambda_1}  \bar{u}_{\lambda_2} u_{\lambda_3} \|_{L^2}
\\ & \lesssim  (\lambda_1 \lambda_2 \lambda_3)^{-s} \lambda^{3s} 
\| u_{\lambda_1}\|_{X^s}\| u_{\lambda_2}\|_{X^s}\| u_{\lambda_3}\|_{X^s}
\end{split}
\]
The summation with respect to the $\lambda_j$'s is straightforward.

{\bf Case 2.} $\max\{\lambda_1,\lambda_2,\lambda_3\}= \mu \gg
\lambda$. In order to have any output at frequency $\lambda$ we must
have at least two $\lambda_j$'s of size $\mu$. Hence we can assume
that
\[
\{ \lambda_1,\lambda_2,\lambda_3\} = \{ \alpha,\mu,\mu\} \qquad 
\alpha \lesssim \mu
\]
We consider two possibilities:

{\bf Case 2a.} $ \alpha \lesssim \lambda \ll \mu$. We begin with the
bound
\begin{equation}
\| P_\lambda(\chi_{[0,1]} v_{\lambda_1}  \bar{v}_{\lambda_2} v_{\lambda_3})\|_{DU^2_\lambda} 
\lesssim \mu^{-1} \log \lambda  \|v_{\lambda_1}\|_{U^2_{\lambda_1}} 
 \|  v_{\lambda_2}\|_{U^2_{\lambda_2}}  \|  v_{\lambda_3}\|_{U^2_{\lambda_3}}
\label{hds}\end{equation}
By duality this is equivalent to 
\[
\left| \int  \chi_{[0,1]} v_{\lambda_1}  \bar{v}_{\lambda_2} v_{\lambda_3}
\bar v_\lambda dxdt \right|   \lesssim  \mu^{-1} \log \lambda  \|  v_{\lambda_1}\|_{U^2_{\lambda_1}} 
 \|  v_{\lambda_2}\|_{U^2_{\lambda_2}}  \|  v_{\lambda_3}\|_{U^2_{\lambda_3}}
\| v_\lambda \|_{V^2_{\lambda}}
\] 
which follows from the bilinear $L^2$ estimate \label{bil2} for the
factors $u_\alpha u_\mu$ and $ \chi_{[0,1]} u_\lambda u_\mu$.

The frequency $\mu$ functions are only controlled on $\mu^{4s}$ time
intervals.  Hence we need to use \eqref{hds} on each such time
interval and then sum up the output from about $\lambda^{4s}
\mu^{-4s}$ such intervals. For each interval $|J| = \lambda^{-4s}$ we
obtain
\[
\lambda^s \| P_\lambda (\chi_J u_{\lambda_1}  \bar{u}_{\lambda_2} u_{\lambda_3})
\|_{DU^2_\lambda}  \lesssim  \alpha^{-s} \mu^{-2s}
\lambda^{5s} \mu^{-4s}  \mu^{-1} \log \lambda\|  u_{\lambda_1}\|_{X^s} 
 \|  u_{\lambda_2}\|_{X^s}  \|  u_{\lambda_3}\|_{X^s}
\]
Then we sum this up with respect to $\alpha$ and $\mu$. This imposes
the restriction $s \geq -\frac16$ but only due to very large values of
$\mu$. We note that we gain almost $1+2s$ derivatives in this
computation.

{\bf Case 2b.} $\alpha \gg \lambda$. For later use we summarize the result 
in this case in the following
\begin{lemma}
\label{lhhh}
Let $I$ be an interval of length $\lambda^{4s}$. Set 
\[
f =   P_\lambda \chi_I
\sum_{\lambda_1,\lambda_2,\lambda_3 \gg \lambda} u_{\lambda_1} \bar
u_{\lambda_2} u_{\lambda_3}
\]
Then we have 
the estimates 
\[
\| Q_{\geq \lambda^2}  f  \|_{X^{0,-\frac12,1}} \lesssim
\lambda^{-1-3s} \|u\|_{X^s}^3
\]
respectively
\[
\|Q_{\leq \lambda^2} f \|_{L^1_t L^2_x+ \lambda^{-\frac14}
L^{\frac43}_x L^1_t} \lesssim \lambda^{-1-3s} \|u\|_{X^s}^3
\]
\end{lemma}

\begin{remark}
 The same estimates remain true, and in fact become easier,  if we replace 
$\overline{u_{\lambda_2}}$ by $u_{\lambda_2}$. 
\end{remark}

We notice that due to the embedding \eqref{embeddingX} and to
Corollary~\ref{lpdual} Lemma \ref{lhhh} implies that
\[
\lambda^s \|f\|_{DU^2_\lambda} \lesssim \lambda^{-1-2s} \|u\|_{X^s}^3
\]
which is a gain similar to the one in Case 2(a). Then the proof of
Proposition~\ref{pnonlin} is concluded.

\begin{proof}[Proof of Lemma~\ref{lhhh}]
To understand the main feature of this case we denote by
$(\tau_i,\xi_i)$ the frequencies for each factor and by $(\tau,\xi)$
the frequency of the output. Then we must have
\[
\xi_1+ \xi_3 = \xi_2+\xi, \qquad \tau_1+ \tau_3 = \tau_2+\tau
\]
This yields 
\[
(\tau_1-\xi_1^2)- (\tau_2-\xi_2^2) + (\tau_3-\xi_3^2) - (\tau-\xi^2) =
2\xi_1\xi_3 - 2 \xi_2 \xi
\]
Since the size of the frequencies $\{ \xi,\xi_1,\xi_2,\xi_3\}$ is $\{
\lambda, \alpha, \mu,\mu\}$ with $\lambda \ll \alpha \lesssim \mu$ we
conclude that
\[
|\tau_1-\xi_1^2|+ |\tau_2-\xi_2^2| + |\tau_3-\xi_3^2| + |\tau-\xi^2|
\gtrsim \alpha \mu, \qquad \text{ if } \lambda_2 = \mu
\]
respectively
\[
|\tau_1-\xi_1^2|+ |\tau_2-\xi_2^2| + |\tau_3-\xi_3^2| + |\tau-\xi^2|
\gtrsim \mu^2, \qquad \text{ if } \lambda_2 = \alpha
\]
This shows that at least one modulation has to be large, namely at
least $\alpha \mu$.  To take advantage of this we split each factor
into a low and a high modulation component.  There are several cases
to consider:

{\bf Case I}. This is when we have three small modulations. Then the output
has large modulation. Depending 
on whether the conjugated factor has lower frequency or not we 
divide this case in three:

{\bf Case I(a)} Here we consider the first component of $f$, namely
\[
f_1 = \sum_{\lambda \ll \alpha \ll \mu} f_1^{\alpha \mu} 
= \sum_{\lambda \ll \alpha \ll \mu} P_\lambda(Q_{\ll
  \alpha \mu} (\chi_I u_{\mu}) \overline
{Q_{\ll \alpha \mu}(\chi_I u_{\mu})} Q_{\ll
  \alpha \mu}(\chi_I u_{\alpha}))
\]
Then $ f_1^{\alpha \mu}$ is localized at modulation $\alpha \mu$.
We begin with an $L^2$ bound for the triple product $P_\lambda
(v_{\mu} \bar v_\mu v_\alpha) $.  We claim that
\begin{equation}
\| P_\lambda (v_{\mu} \bar v_\mu v_\alpha) \|_{L^2} \lesssim
\lambda^{\frac12} \mu^{-\frac12} \| v_{\mu}\|_{U^2_\mu} \| v_\mu\|_{U^2_\mu}  \| v_\alpha \|_{U^2_\alpha}
\label{tril2}\end{equation}
Indeed, using the energy bound for $v_{\mu}$ and the 
bilinear $L^2$ bound for $\bar v_\mu v_\alpha$ we obtain
\[
\|v_{\mu} \bar v_\mu v_\alpha \|_{L^2 L^1} \lesssim \mu^{-\frac12} \|
v_{\mu}\|_{U^2_\mu} \| v_\mu\|_{U^2_\mu} \| v_\alpha \|_{U^2_\alpha}
\]
Applying $P_\lambda$ the estimate \eqref{tril2}  follows from
Bernstein's inequality.

To use \eqref{tril2} in order to bound $f_1$ we decompose each factor
$u_\mu$, respectively $u_\lambda$ with respect to time intervals of
length $\mu^{4s}$, respectively $\alpha^{4s}$ and apply \eqref{tril2}
for each combination. The contributions of $\mu^{4s}$ separated
intervals is negligible since the kernel of $Q_{\ll \alpha \mu}$
decays rapidly on the $(\alpha \mu)^{-1}$ timescale. Hence there
are about $\lambda^{4s} \mu^{-4s}$ contributions to add up.
We obtain 
\[
\| f^{\mu\alpha}_{1} \|_{L^2} \lesssim \lambda^{4s}
\mu^{-4s} \alpha^{-s} \mu^{-2s} \lambda^{\frac12} \mu^{-\frac12}
\| u_{\mu}\|_{X^s} \| u_\mu\|_{X^s}  \| u_\alpha \|_{X^s}
\]
Since $ f^{\mu\alpha}_{1}$ has modulation $\alpha \mu$ this gives
\[
\| f^{\mu\alpha}_{1} \|_{X^{0,-\frac12,1}} \lesssim \lambda^{\frac12+4s}
\mu^{-1-6s} \alpha^{-\frac12-s} 
\| u_{\mu}\|_{X^s} \| u_\mu\|_{X^s}  \| u_\alpha \|_{X^s}
\]
The summation with respect to the dyadic indices $\alpha$ and $\mu$ is
straightforward provided that $s \geq -\frac16$. We obtain
\[
\| f_{1} \|_{X^{0,-\frac12,1}} \lesssim \lambda^{-1-3s} \|u\|_{X^s}^3
\]

{\bf Case I(b)} The second component of $f$ is
\[
f_2 = \sum_{\lambda \ll \alpha \ll \mu} f_2^{\alpha \mu} 
= \sum_{\lambda \ll \alpha \ll \mu} P_\lambda(Q_{\ll
  \mu^2} (\chi_I u_{\mu}) \overline
{Q_{\ll  \mu^2}(\chi_I u_{\alpha})} Q_{\ll
   \mu^2}(\chi_I u_{\mu}))
\]
Then $ f_2^{\alpha \mu}$ is localized at modulation $\mu^2$.  We can
still use \eqref{tril2} since the location of the complex conjugates
does not matter. Hence $ f_2^{\alpha \mu}$ satisfies the same $L^2$
bound as $f_1^{\alpha \mu}$.  However, because of the larger
modulation we obtain a better $X^{0,-\frac12,1}$ bound, namely
\[
\| f^{\mu\alpha}_{2} \|_{X^{0,-\frac12,1}} \lesssim \lambda^{\frac12+4s}
\mu^{-\frac32-6s} \alpha^{-s} 
\| u_{\mu}\|_{X^s} \| u_\mu\|_{X^s}  \| u_\alpha \|_{X^s}
\]
After summation with respect to $\alpha$ and $\mu$ we obtain the same
bound for $f_2$ as for $f_1$; the difference is that the summation can
be carried out for $s \geq -\frac3{14}$.

{\bf Case I(c)} The third component of $f$ is
\[
f_3 = \sum_{\lambda \ll  \mu} f_3^{\mu} 
= \sum_{\lambda \ll \mu} P_\lambda(Q_{\ll
  \mu^2} (\chi_I u_{\mu}) \overline
{Q_{\ll  \mu^2}(\chi_I u_{\mu})} Q_{\ll
   \mu^2}(\chi_I u_{\mu}))
\]
Then $ f_3^{\alpha \mu}$ is localized at modulation $\mu^2$.  We claim
that \eqref{tril2} still holds.  To prove this we first observe that
in order for the output to be at low frequency $\lambda$, two of the
frequencies $\xi_1,-\xi_2,\xi_3$ must be $\mu$ separated. Then we  use the
bilinear $L^2$ bound for those two factors, and the energy bound for
the third.

By \eqref{tril2} we obtain as in Case I(a)
\[
\| f_3^{\mu}\|_{X^{0,-\frac12,1}} \lesssim \lambda^{\frac12+4s}
\mu^{-\frac3{2}-7s} \| u_{\mu}\|_{X^s} \| u_\mu\|_{X^s} \| u_\mu
\|_{X^s}
\]
and the summation with respect to $\mu$ can be carried out for $s \geq
-\frac3{14}$.

{\bf Case II}. This is when at least one factor has large modulation.
Depending on which factor has large modulation and on whether the
conjugated factor has lower frequency or not we divide this case in
six:

{\bf Case II(a)}. Here we consider 
\[
f_4 = \sum_{\lambda \ll \alpha \ll \mu} f_4^{\alpha \mu} 
= \sum_{\lambda \ll \alpha \ll \mu} P_\lambda (\chi_I  Q_{\gtrsim
  \alpha \mu} (\chi_I u_{\mu}) \overline
{   u_{\mu}}   u_{\alpha})  +  P_\lambda ( \chi_I   u_{\mu} \overline
{Q_{\gtrsim \alpha \mu}(\chi_I u_{\mu})}   u_{\alpha})
\]
The two terms are similar, so we restrict our attention to the first
one. Our starting point is the bound 
\begin{equation}
\| Q_{\gtrsim
  \alpha \mu} v_\mu \overline
{v_{\mu}} v_{\alpha}\|_{L^1} \lesssim \alpha^{-\frac12}  \mu^{-1}
\|v_\mu\|_{U^2_\mu}
 \|{v_{\mu}}\|_{U^2_\mu} \| v_{\alpha}\|_{U^2_\alpha}
\label{hll}\end{equation}
which is obtained from the $L^2$ estimate for the first factor and a
bilinear $L^2$ estimate for the remaining product. 

{\bf Low modulation output:}
By Bernstein's inequality \eqref{hll} implies
\begin{equation}
  \| Q_{\gtrsim
    \alpha \mu} v_\mu \overline
  {v_{\mu}} v_{\alpha}\|_{L^1L^2} \lesssim \lambda^{\frac12} 
\alpha^{-\frac12}  \mu^{-1}
  \|v_\mu\|_{U^2_\mu}
  \|{v_{\mu}}\|_{U^2_\mu} \| v_{\alpha}\|_{U^2_\alpha}
\label{hlla}\end{equation}
 Summing up \eqref{hlla} over $\lambda^{4s}
\mu^{-4s}$ time intervals of length $\mu^{4s}$ we obtain
\[
\|f_{4}^{\alpha \mu}\|_{L^1L^2} \lesssim \alpha^{-s} \mu^{-2s} 
\lambda^{4s} \mu^{-4s}  \lambda^\frac12 \alpha^{-\frac12}  \mu^{-1} 
\|u_\mu\|_{X^s}
 \|{u_{\mu}}\|_{X^s} \| u_{\alpha}\|_{X^s}
\]
For $s \geq -\frac16$ we can sum this up with respect to $\alpha$
and $\lambda$ to obtain
\[
\|f_{4}\|_{L^1 L^2}  \lesssim \lambda^{-1 - 3s} \|u\|_{X^s}^3
\]

{\bf Intermediate modulation output:} Consider now the
$X^{0,-\frac12,1}$ estimate at modulations $\lambda^2 \leq \sigma \leq
\alpha \mu$. From \eqref{hll} and Bernstein's inequality we obtain
\begin{equation}
\| Q_\sigma P_\lambda (Q_{\gtrsim
  \alpha \mu} v_\mu \overline
{v_{\mu}} v_{\alpha})\|_{L^2} \lesssim (\lambda \sigma)^\frac12
 \alpha^{-\frac12}  \mu^{-1}
\|v_\mu\|_{U^2_\mu}
 \|{v_{\mu}}\|_{U^2_\mu} \| v_{\alpha}\|_{U^2_\alpha}
\label{hll1}\end{equation}
The kernel of $Q_\sigma$ is rapidly decaying off diagonal on the
$\sigma^{-1}$ scale. Then in estimating the sum over $\mu^{4s}$
intervals there is a gain coming from the fact that we only need
square summability with respect to intervals of size $\sigma^{-1}$. We
consider two cases.

a) If $\sigma^{-1} < \mu^{4s}$ then we need square summability with
respect to intervals of size $\mu^{4s}$ so we obtain
\[
\| Q_\sigma f_{4}^{\alpha \mu} \|_{L^2} \lesssim \alpha^{-s}
\mu^{-2s} \lambda^{2s} \mu^{-2s} (\lambda \sigma)^\frac12
\alpha^{-\frac12} \mu^{-1} \|u_\mu\|_{X^s} \|{u_{\mu}}\|_{X^s} \|
u_{\alpha}\|_{X^s}
\]
or equivalently
\[
\| Q_\sigma f_{4}^{\alpha \mu} \|_{X^{0,-\frac12,1}} \lesssim 
\lambda^{\frac12+2s} \alpha^{-\frac12-s} \mu^{-1-4s} \|u_\mu\|_{X^s}
\|{u_{\mu}}\|_{X^s} \| u_{\alpha}\|_{X^s}
\]
Adding up with respect to $\sigma$ yields
\[ 
\sum_{\sigma = \max\{\lambda^2,\mu^{-4s}\}}^{\alpha \mu}
 \| Q_\sigma f_{4}^{\alpha \mu} \|_{X^{0,-\frac12,1}}  \lesssim   
\lambda^{\frac12+2s} \alpha^{-\frac12-s} \mu^{-1-4s}
\ln \mu
  \|u_\mu\|_{X^s}
\|{u_{\mu}}\|_{X^s} \| u_{\alpha}\|_{X^s}
\]
and now the summation with respect to $\alpha$ and $\mu$ is
straightforward for $s > -\frac14$.

b) If $\sigma^{-1} > \mu^{4s}$ then we need square summability with
respect to intervals of size $\sigma^{-1}$
so we obtain
\[
\| Q_\sigma f_{21}^{\alpha \mu} \|_{L^2} \lesssim \alpha^{-s}
\mu^{-2s} \lambda^{2s} \sigma^{-\frac12} \mu^{-4s} (\lambda
\sigma)^\frac12 \alpha^{-\frac12} \mu^{-1} \|u_\mu\|_{X^s}
\|{u_{\mu}}\|_{X^s} \| u_{\alpha}\|_{X^s}
\]
or equivalently
\[
\| Q_\sigma f_{21}^{\alpha \mu} \|_{X^{0,-\frac12,1}} \lesssim
\sigma^{- \frac12} \lambda^{\frac12+2s} \alpha^{-\frac12-s} \mu^{-1-6s}
\|u_\mu\|_{X^s} \|{u_{\mu}}\|_{X^s} \| u_{\alpha}\|_{X^s}
\]
Adding up  with respect to $\sigma$ gives
\[
\sum_{\sigma = \lambda^2}^{\mu^{-4s}} \| Q_{\sigma} 
 f_{4}^{\alpha \mu} \|_{X^{0,-\frac12,1}} \lesssim 
\lambda^{-\frac12+2s} \alpha^{-\frac12-s} \mu^{-1-6s}
\|u_\mu\|_{X^s}
\|{u_{\mu}}\|_{X^s} \| u_{\alpha}\|_{X^s}
\]
In this case the summation with respect to $\alpha, \mu$ 
gains $-2+4s$ derivatives,
which is
better result than needed,
 but the summation requires $s \geq -\frac16$.

{\bf High modulation output:} Here we estimate the output at
modulations $\sigma \gg \alpha \mu$.  In order to obtain such an
output at least one of the factors must have modulation at least
$\sigma$. Without any restriction in generality we assume that this is
the first factor, as the other cases are considerably simpler. 
This gives a trivial improvement in the $L^2$ bound for
the first factor, so instead of \eqref{hll1} we have
\begin{equation}
\| Q_\sigma P_\lambda (Q_{\gtrsim
  \alpha \mu} v_\mu \overline
{v_{\mu}} v_{\alpha})\|_{L^2} \lesssim (\lambda \sigma)^\frac12
 \sigma^{-\frac12}  \mu^{-\frac12}
\|v_\mu\|_{U^2_\mu}
 \|{v_{\mu}}\|_{U^2_\mu} \| v_{\alpha}\|_{U^2_\alpha}
\label{hll1a}\end{equation}
Then the bound in case (a) above is replaced by
\[
\| Q_\sigma f_{4}^{\alpha \mu} \|_{L^2} \lesssim \alpha^{-s}
\mu^{-2s} \lambda^{2s} \mu^{-2s} (\lambda \sigma)^\frac12
\sigma^{-\frac12} \mu^{-\frac12} \|u_\mu\|_{X^s} \|{u_{\mu}}\|_{X^s} \|
u_{\alpha}\|_{X^s}
\]
or equivalently
\[
\| Q_\sigma f_{21}^{\alpha \mu} \|_{X^{0,-\frac12,1}} \lesssim 
\lambda^{\frac12+2s} \alpha^{-s}  \mu^{-\frac12-4s} \sigma^{-\frac12} \|u_\mu\|_{X^s}
\|{u_{\mu}}\|_{X^s} \| u_{\alpha}\|_{X^s}
\]
which has better summability with respect to large $\sigma$.

{\bf Case II(b).} This is when the low frequency factor has high
modulation. We consider terms of the form
\[
f_5 = \sum_{\lambda \ll \alpha \ll \mu} f_5^{\alpha \mu} =  
\sum_{\lambda \ll \alpha \ll \mu}
P_\lambda ( \chi_I    u_{\mu} \overline
{ u_{\mu}} Q_{\gtrsim
  \alpha \mu}(\chi_I u_{\alpha}))
\]
Depending on the relative size of $\alpha$ and $\mu$ we 
divide the problem into two sub-cases:

{\bf Case II(b)-1.} $\lambda \mu \leq \alpha^2$.  By orthogonality we
can assume that the two $u_\mu$ factors are frequency localized in
$\alpha$ separated intervals of length $\alpha$.  Then we use the
bilinear $L^2$ bound for their product and the $L^2$ bound for the
high modulation factor to obtain a weaker analogue of \eqref{hll},
namely
\begin{equation}
  \| P_\lambda(v_\mu \overline
  {v_{\mu}}  Q_{\gtrsim
    \alpha \mu} v_{\alpha})\|_{L^1} \lesssim \alpha^{-1}  \mu^{-\frac12}
  \|v_\mu\|_{U^2_\mu}
  \|{v_{\mu}}\|_{U^2_\mu} \| v_{\alpha}\|_{U^2_\alpha}
\label{llh}\end{equation}

{\bf Low modulation output.} Compensating for the weaker bound
\eqref{llh}, in this case there is an improvement in the summation
over time intervals.  We decompose the $\lambda^{4s}$ time interval
$I$ in two steps.  First we split it into $\lambda^{4s} \alpha^{-4s}$
time intervals of length $ \alpha^{4s}$, which gives a $\lambda^{4s}
\alpha^{-4s}$ factor in the summation. Secondly we split each
$\alpha^{4s}$ time interval into $\alpha^{4s} \mu^{-4s}$ time
intervals of length $ \mu^{4s}$. Since the frequency $\alpha$ factor
is square summable with respect to this partition, by Cauchy-Schwartz
this gives only an $\alpha^{2s} \mu^{-2s}$ factor in the summation. We
obtain
\[
\|f_{5}^{\alpha \mu}\|_{L^1} \lesssim \alpha^{-s} \mu^{-2s} 
\lambda^{4s} \alpha^{-2s} \mu^{-2s}  \alpha^{-1}  \mu^{-\frac12} 
\|u_\mu\|_{X^s}
 \|{u_{\mu}}\|_{X^s} \| u_{\alpha}\|_{X^s}
\]
which by Bernstein's inequality implies that
\[
\| Q_{\leq \lambda^2} f_{5}^{\alpha \mu}\|_{L^1 L^2} \lesssim 
\lambda^{\frac12+ 4s} \alpha^{-1-3s}  \mu^{-\frac12-4s} 
\|u_\mu\|_{X^s}
 \|{u_{\mu}}\|_{X^s} \| u_{\alpha}\|_{X^s}
\]
This is summable with respect to large $\mu$ only if $s \geq
-\frac18$. However, the restriction $\lambda \mu \leq \alpha^2$
improves the summation. Assuming $s < -\frac18$ the $\mu$ summation
yields
\[
\sum_{\alpha \leq \mu \leq \lambda^{-1} \alpha^2} \| Q_{\leq
  \lambda^2} f_{5}^{\alpha \mu}\|_{L^1 L^2} \lesssim \lambda^{1+
  8s} \alpha^{-2-11s}  \|u\|_{X^s}^2  \| u_{\alpha}\|_{X^s}
\]
which is summable with respect to $\alpha$ for $s > -\frac2{11}$.

{\bf Intermediate modulation output.} $\lambda^2 \leq \sigma \leq
\alpha \mu$. Here we argue as in Case II(a) but using \eqref{llh}
instead of \eqref{hll}.  From \eqref{llh} and Bernstein's inequality we obtain
\begin{equation}
\| Q_\sigma P_\lambda ( v_\mu \overline
{v_{\mu}} Q_{\gtrsim
  \alpha \mu} v_{\alpha})\|_{L^2} \lesssim (\lambda \sigma)^\frac12
 \alpha^{-1}  \mu^{-\frac12}
\|v_\mu\|_{U^2_\mu}
 \|{v_{\mu}}\|_{U^2_\mu} \| v_{\alpha}\|_{U^2_\alpha}
\label{lhh0}\end{equation}
We split this again depending on $\sigma$ but also taking into account
the improved summability up to the $\alpha^{4s}$ time scale, as
discussed above for the case of  low modulation output.

a) If $\sigma^{-1} \leq \mu^{4s}$ then due to the square integrability
of $u_\alpha$ in each $\lambda^{4\alpha}$ time interval we have an
interval summation factor $\lambda^{2s} \alpha^{-2s}$.  Hence
\[
\| Q_\sigma f_5^{\alpha \mu} \|_{L^2} \lesssim\alpha^{-s}
\mu^{-2s}   \lambda^{2s}  \alpha^{-2s}     
 (\lambda \sigma)^\frac12
 \alpha^{-1}  \mu^{-\frac12} 
\|u_\mu\|_{X^s}
 \|{u_{\mu}}\|_{X^s} \| u_{\alpha}\|_{X^s}
\]
which yields
\[
\| Q_\sigma f_5^{\alpha \mu} \|_{X^{0,-\frac12,1}} \lesssim 
\lambda^{\frac12+2s} \alpha^{-1-3s} \mu^{-\frac12-2s}
\|u_\mu\|_{X^s} \|{u_{\mu}}\|_{X^s} \|
u_{\alpha}\|_{X^s}
\]
The summation with respect to $\sigma$, $\alpha$ and $\mu$ is
straightforward for $s \geq -\frac14$.

b) The case $\lambda^2 < \sigma < \mu^{-4s}$ is somewhat worse because
the kernel of $Q_\sigma$ decays only on the $\sigma^{-1}$  scale 
which is now larger than $\mu^{4s}$. Hence inputs from $\mu^{4s}$ time
intervals within each $\sigma^{-1}$ time interval are no longer orthogonal.
This yields a weaker interval summation factor, namely $\lambda^{2s}
\alpha^{-2s} \sigma^{-\frac12} \mu^{2s}$.  Hence
\[
\| Q_\sigma f_5^{\alpha \mu} \|_{L^2} \lesssim \alpha^{-s}
\mu^{-2s}   \lambda^{2s}  \alpha^{-2s} \sigma^{-\frac12}  \mu^{-2s}     
 (\lambda \sigma)^\frac12
 \alpha^{-1}  \mu^{-\frac12} 
\|u_\mu\|_{X^s}
 \|{u_{\mu}}\|_{X^s} \| u_{\alpha}\|_{X^s}
\]
which yields
\[
\| Q_\sigma f_5^{\alpha \mu} \|_{X^{0,-\frac12,1}} 
\lesssim \sigma^{-\frac12}
\lambda^{\frac12+2s} \alpha^{-1-3s} \mu^{-\frac12-4s}
\|u_\mu\|_{X^s} \|{u_{\mu}}\|_{X^s} \|
u_{\alpha}\|_{X^s}
\]

The summation with respect to $\sigma$ is straightforward:
\[
\sum_{\sigma= \lambda^2}^{\mu^{-4s}}  
\| Q_\sigma f_5^{\alpha \mu} \|_{X^{0,-\frac12,1}} \lesssim 
\lambda^{-\frac12+2s} \alpha^{-1-3s} \mu^{-\frac12-4s}
\|u_\mu\|_{X^s} \|{u_{\mu}}\|_{X^s} \|u_{\alpha}\|_{X^s}
\]
 However, in
the $\alpha$ and $\mu$ summation we need to use the restriction
$\lambda \mu \leq \alpha^2$ exactly as in the case of low modulation
output.

{\bf High modulation output:} Here we estimate the output at
modulations $\sigma \gg \alpha \mu$. Then we can assume that the last
factor has modulation at least $\sigma$ therefore it satisfies a
better $L^2$ bound, which leads to
\begin{equation}
\| Q_\sigma P_\lambda ( v_\mu \overline
{v_{\mu}} Q_{\gtrsim
  \alpha \mu} v_{\alpha})\|_{L^2} \lesssim \lambda^\frac12
 \alpha^{-\frac12} 
\|v_\mu\|_{U^2_\mu}
 \|{v_{\mu}}\|_{U^2_\mu} \| v_{\alpha}\|_{U^2_\alpha}
\label{lhhah}\end{equation}
Then instead of the
estimate in case (a) above we obtain
\[
\| Q_\sigma f_{5}^{\alpha \mu} \|_{L^2} \lesssim \alpha^{-s}
\mu^{-2s} \lambda^{2s} \alpha^{-2s} \lambda^\frac12
 \alpha^{-\frac12} \|u_\mu\|_{X^s} \|{u_{\mu}}\|_{X^s} \|
u_{\alpha}\|_{X^s}
\]
or equivalently
\[
\| Q_\sigma f_{5}^{\alpha \mu} \|_{X^{0,-\frac12,1}} \lesssim
\lambda^{\frac12+2s} \alpha^{-\frac12-3s} \mu^{-2s} \sigma^{-\frac12}
\|u_\mu\|_{X^s} \|{u_{\mu}}\|_{X^s} \| u_{\alpha}\|_{X^s}
\]
and hence 
\[
\| Q_{\gg \alpha \mu} f_5^{\alpha \mu} \|_{X^{0,-\frac12,1}} 
\lesssim \lambda^{-1/2 +2s} \alpha^{-1-3s} \mu^{-\frac12 -4s} 
 \|u_\mu\|_{X^s}
\|{u_{\mu}}\|_{X^s} \| u_{\alpha}\|_{X^s}.
\]
 The
condition $\lambda \mu \leq \alpha^2$ is again needed.

{\bf Case II(b)-2:}  $\lambda \mu > \alpha^2$.  Then the arguments in
the previous case fail to provide enough decay in order to insure
summability for very large $\mu$.

{\bf Low modulation output}.
In this case we are able to establish the following improvement of \eqref{llh},
\begin{equation}
\| Q_{< \lambda^2} P_\lambda(v_\mu \overline
  {v_{\mu}}  Q_{\gtrsim
    \alpha \mu} v_{\alpha})\|_{L^1 L^2} \lesssim  \mu^{-1}
  \|v_\mu\|_{U^2_\mu}
  \|{v_{\mu}}\|_{U^2_\mu} \| v_{\alpha}\|_{U^2_\alpha}
\label{llh1}\end{equation}
The rest of the analysis is similar to the computation in Case
II(b)-1.  The only difference is that here we gain an extra factor of 
$\alpha (\lambda \mu)^{-\frac12} \leq 1$, which improves the summation
for large $\mu$.

To prove \eqref{llh1} we only use the $L^2$ bound for $v_{\alpha}$.
Then, using the atomic decomposition for each of the two $v_\mu$
factors, we conclude that it suffices to prove \eqref{llh1} in the
case when both $v_\mu$ factors solve the linear equation.  By
orthogonality we can assume that both are frequency localized in
$\alpha$ intervals which are $\alpha$ separated.  Then we use the
$L^2$ bound for the product of $u_\mu \bar u_\mu$,
\[
\|v_\mu \bar v_\mu\|_{L^2} \lesssim \alpha^{-\frac12}  \|v_\mu\|_{U^2_\mu}
  \|{v_{\mu}}\|_{U^2_\mu}
\]

However, due to the frequency localization we also obtain that the
product is Fourier localized in a thin rectangle $R$ of size
$\alpha^2/\mu \times \alpha \mu$ at slope $\mu^{-1}$. Next we
consider the product 
\[
(v_\mu \bar v_\mu) \cdot (Q_{\gtrsim \alpha \mu} v_{\alpha})
\]
which we view as a product of two $L^2$ functions with different
Fourier localizations. The product is only estimated in a Fourier
rectangle of size $\lambda \times \lambda^2$, therefore by
orthogonality it suffices to estimate the product assuming that both
factors are Fourier localized in rectangles $R_1$, $R_2$ of similar
size.  The intersection $R_0= R \cap R_1$ is a shorter
rectangle  of size $\alpha^2/\mu \times \lambda^2$. Our
assumption $\alpha^2 < \lambda \mu$ insures that $R_0$ is essentially
vertical.  But by Bernstein's inequality 
we have the pointwise bound
\[
\| g\|_{L^2 L^\infty} \lesssim \alpha \mu^{-\frac12} \|g\|_{L^2},
\qquad \text{supp}\  \hat g \subset R_0
\]
therefore \eqref{llh1} follows.

{\bf Intermediate modulation output,} $\lambda^2 < \sigma \lesssim
\alpha \mu$. Then a similar argument applies. $R_0$ has size
$\alpha^2/\mu \times \sigma$, which yields the pointwise bound
\[
\| g\|_{L^\infty} \lesssim \alpha \mu^{-\frac12} \sigma^\frac12
\|g\|_{L^2}, \qquad  \text{supp}\  \hat g \subset R_0
\]
This in turn leads to
\begin{equation}
  \| Q_{\sigma} P_\lambda(v_\mu \overline
  {v_{\mu}}  Q_{\gtrsim
    \alpha \mu} v_{\alpha})\|_{L^2} \lesssim \sigma^\frac12 
\alpha^{\frac12}  \mu^{-1}
  \|v_\mu\|_{U^2_\mu}
  \|{v_{\mu}}\|_{U^2_\mu} \| v_{\alpha}\|_{U^2_\alpha}
\label{llh3}\end{equation}
which is again  an improvement of 
$\alpha (\lambda \mu)^{-\frac12}$ over the similar computation in Case II(b)-1. 

{\bf Large modulation output,} $ \sigma \gg \alpha \mu$. Then we can
assume that the third factor has modulation at least $\sigma$. Also
$R$ has size $\alpha^2/\mu \times \alpha \mu$, therefore
\[
\| g\|_{L^\infty} \lesssim \alpha^\frac32 
\|g\|_{L^2}, \qquad  \text{supp}\  \hat g \subset R
\]
which implies that
\begin{equation}
  \| Q_{\sigma} P_\lambda(v_\mu \overline
  {v_{\mu}}  Q_{\gtrsim
    \sigma} v_{\alpha})\|_{L^2} \lesssim \sigma^{-\frac12} 
\alpha 
  \|v_\mu\|_{U^2_\mu}
  \|{v_{\mu}}\|_{U^2_\mu} \| v_{\alpha}\|_{U^2_\alpha},
\label{llh4}\end{equation}
an improvement of at least $\alpha (\lambda \mu)^{-\frac12}$ over 
\eqref{lhhah}. The conclusion follows in a similar fashion.

 {\bf Case II(c).} This is when the low frequency factor is conjugated
but does not have high modulation. We consider terms of the form
\[
f_6 = \sum_{\lambda \ll \alpha \lesssim \mu} f_{6}^{\alpha \mu}
= \sum_{\lambda \ll \alpha \lesssim \mu} P_\lambda (\chi_I 
Q_{\gtrsim
  \mu^2} (\chi_I u_{\mu}) \overline
{  u_{\alpha}}   u_{\mu})
\]
If $\alpha \ll \mu$ then the last two factors are $\mu$ separated in
frequency. But even if $\alpha \approx \mu$, in order for the final
output to be at frequency $\lambda$ the  two last factors must still
be $\mu$ separated. Then we can use the trilinear bound
\begin{equation}
\| P_\lambda(Q_{\gtrsim
  \alpha \mu} v_\mu \overline
{v_{\alpha}} v_{\mu})\|_{L^1} \lesssim   \alpha^{-1/2} \mu^{-1}
\|v_\mu\|_{U^2_\mu}
 \|{v_{\mu}}\|_{U^2_\mu} \| v_{\alpha}\|_{U^2_\alpha},
\label{lhl}\end{equation} 
obtained by estimating in $L^2$ the first factor and the
remaining product.

 The constants here are better than the 
ones in Case II(a), and the rest of the argument proceeds
as there without any significant changes.

 {\bf Case II(d).} This is when the low frequency factor is conjugated
and  has high modulation. We consider terms of the form
\[
f_7 = \sum_{\lambda \ll \alpha \ll \mu} f_{7}^{\alpha \mu}
= \sum_{\lambda \ll \alpha \ll \mu} P_\lambda (\chi_I 
 u_{\mu} \overline
{  Q_{\gtrsim
  \mu^2} (\chi_I u_{\alpha})}   u_{\mu})
\]
In order for the final output to be at frequency $\lambda$ the two
frequency $\mu$ factors must still be $\mu$ separated. This leads to
the trilinear bound
\begin{equation}
\| P_\lambda( v_\mu \overline
{Q_{\gtrsim
  \alpha \mu} v_{\alpha}} v_{\mu})\|_{L^1} \lesssim   \mu^{-\frac32}
\|v_\mu\|_{U^2_\mu}
 \|{v_{\mu}}\|_{U^2_\mu} \| v_{\alpha}\|_{U^2_\alpha},
\label{lhla}\end{equation} 
and the argument is completed again as in Case II(a) but with better
constants.

{\bf Case II(e).} This is when all frequencies are equal and the
conjugated  factor has high modulation.
 We consider terms of the form
\[
f_8 = \sum_{\lambda  \ll \mu} f_{8}^{\mu}
= \sum_{\lambda  \ll \mu} P_\lambda (\chi_I 
 u_{\mu} \overline
{  Q_{\gtrsim
  \mu^2} (\chi_I u_{\mu})}   u_{\mu})
\] 
In some sense this is the worst case because we cannot enforce any
frequency separation between the two unconjugated factors.  We still
want to gain some power of $\mu$ in order to have summability for
large $\mu$.  At least to some extent we can do this by the
lateral Strichartz estimates in Corollary~\ref{upstr}(bc) to obtain
\begin{equation}
 \| v_\mu \overline
  { Q_{\gtrsim
     \mu^2 } v_{\mu}} v_{\mu}\|_{L^\frac43_x L^1_t} \lesssim   \mu^{-\frac54}
  \|v_\mu\|_{U^2_\mu}
  \|{v_{\mu}}\|_{U^2_\mu} \| v_{\mu}\|_{U^2_\mu}
\label{hhh}\end{equation}
This is done for instance by using the $L^\infty_x L^2_t$ bound for
one $v_\mu$ factor, respectively the $L^4_x L^\infty_t$ for the other  
$v_\mu$ factor.

{\bf Low modulation output:}
After summation with respect to $\mu^{4s}$ time intervals \eqref{hhh} gives
\[
\|Q_{< \lambda^2} P_\lambda(u_\mu \overline
  { Q_{\gtrsim
    \mu^2} u_{\mu}}  u_{\mu})\|_{L^\frac43_x L^1_t} \lesssim
\lambda^{4s} \mu^{-4s} \mu^{-3s} \mu^{-\frac54} \|u_\mu\|_{X^s}^3
\]
which is easily summed up with respect to $\mu$ for $s \geq
-\frac5{28}$.

{\bf Intermediate modulation output:} $\lambda^2 < \sigma \leq \mu^2$.
>From \eqref{hhh} combined with Bernstein's inequality we
obtain
\begin{equation}
 \| Q_\sigma P_\lambda(v_\mu \overline
  { Q_{\gtrsim
     \mu^2 } v_{\mu}} v_{\mu})\|_{L^2} \lesssim
 \sigma^\frac12 \lambda^\frac14  \mu^{-\frac54}
  \|v_\mu\|_{U^2_\mu}
  \|{v_{\mu}}\|_{U^2_\mu} \| v_{\mu}\|_{U^2_\mu}
\label{hhha}\end{equation}
Adding this up with respect to $\mu^{4s}$ time intervals yields
\[
\|Q_{\sigma} P_\lambda(u_\mu \overline { Q_{\gtrsim \mu^2} u_{\mu}}
u_{\mu})\|_{L^2} \lesssim \lambda^{\frac14+4s} \sigma^\frac12
\mu^{-\frac54 -7s} \|u_\mu\|_{X^s}^3
\]
or equivalently
\[
\|Q_{\sigma} P_\lambda(u_\mu \overline { Q_{\gtrsim \mu^2} u_{\mu}}
u_{\mu})\|_{X^{0,-\frac12,1}} \lesssim \lambda^{\frac14+4s}
\mu^{-\frac54 -7s} \|u_\mu\|_{X^s}^3
\]
which is easily summed up with respect to $\sigma$ and $\mu$  for $s >
-\frac5{28}$.

Finally, the output modulations which are larger than $\mu^2$ are
treated as in the first case.

{\bf High modulation output:} $\sigma \gg \mu^2$. Without any
restriction in generality we assume that the second factor has
modulation at least $\sigma$. Instead of \eqref{hhha} we get
\begin{equation}
 \| Q_\sigma P_\lambda(v_\mu \overline
  { Q_{\gtrsim
     \sigma } v_{\mu}} v_{\mu})\|_{L^2} \lesssim
  \lambda^\frac14  \mu^{-\frac14}
  \|v_\mu\|_{U^2_\mu}
  \|{v_{\mu}}\|_{U^2_\mu} \| v_{\mu}\|_{U^2_\mu}
\label{hhhb}\end{equation}
Adding this up with respect to $\lambda^{4s}\mu^{-4s}$ time intervals yields
\[
\|Q_{\sigma} P_\lambda(u_\mu \overline { Q_{\gtrsim \sigma} u_{\mu}}
u_{\mu})\|_{L^2} \lesssim \lambda^{\frac14+4s} 
\mu^{-\frac14 -7s} \|u_\mu\|_{X^s}^3
\]
or equivalently
\[
\|Q_{\sigma} P_\lambda(u_\mu \overline { Q_{\gtrsim \mu^2} u_{\mu}}
u_{\mu})\|_{X^{0,-\frac12,1}} \lesssim \lambda^{\frac14+4s}
\sigma^{-\frac12} \mu^{-\frac14 -7s} \|u_\mu\|_{X^s}^3
\]
The summation with respect to $\sigma$ and $\mu$ requires again $s >
-\frac5{28}$.

\end{proof}

\section{The energy conservation}

 It remains to study the conservation of the $H^s$ energy.
We first set 
\[
E_0(u) = \langle A(D) u,u\rangle
\]
For the straight $H^s$ energy conservation it suffices to take
\[
 a(\xi) = (1+\xi^2)^{s} 
\]
However in order to gain the uniformity in $t$ required by \eqref{xs}
we need to allow a slightly larger class of symbols.

\begin{definition}
  Let $s \in \R$ and $\e > 0$. Then $S_\epsilon^s$ is the class of
  spherically symmetric symbols with the following properties:

(i) symbol regularity,
\[
| \partial^\alpha a(\xi)| \lesssim a(\xi) (1+\xi^2)^{-\alpha/2}
\]

(ii) decay at infinity,
\[
s \leq \frac{\ln a(\xi)}{\ln(1+\xi^2)} \leq s+\epsilon, \qquad
s-\epsilon \leq \frac{d \ln a(\xi)} {d \ln(1+\xi^2)} \leq s+\epsilon
\]
\end{definition}
Here $\epsilon$ is a small parameter.

We compute the derivative of $E_0$ along the flow,
\[
\frac{d}{dt} E_0(u) =  R_4(u) = 2 \Re \langle i A(D) u,|u|^2 u\rangle 
\]
We write $R_4$ as a multilinear operator in the Fourier space,
\[
R_4(u) = 2 \Re \int_{P_4} i a(\xi_1) \hat u(\xi_1)
\hat u(\xi_2) \overline{ \hat u(\xi_3)
\hat u(\xi_4)} d\sigma  
\]
where 
\[
P_4 = \{ \xi_1+\xi_2-\xi_3-\xi_4=0\}
\]
This can be symmetrized,
\[
R_4(u) = \frac12 \Re \int_{P_4} i
(a(\xi_1)+a(\xi_2)-a(\xi_3)-a(\xi_4)) \hat u(\xi_1) \hat u(\xi_2) \overline{\hat
u(\xi_3) \hat u(\xi_4)} d\sigma
\]
Following a variation of the $I$-method, see Tao~\cite{MR2233925}-3.9
and references therein, we seek to cancel this term by
perturbing the energy, namely by
\[
E_1(u) =  \int_{P_4}
b_4(\xi_1,\xi_2,\xi_3,\xi_4) \hat u(\xi_1) \hat u(\xi_2) \overline{\hat
u(\xi_3) \hat u(\xi_4)} d\sigma
\]
To determine the best choice for $B$  we compute
\[
\begin{split}
\frac{d}{dt} E_1(u) &=  \int_{P_4} i 
b_4(\xi_1,\xi_2,\xi_3,\xi_4)
(\xi_1^2+\xi_2^2- \xi_3^2- \xi_4^2) \hat u(\xi_1) \hat u(\xi_2) \overline{\hat
u(\xi_3) \hat u(\xi_4)} d\sigma \\ & + R_6(u)
\end{split}
\]
where $R_6(u)$ is given by
\[
R_6(u) = 4 \Re \int _{\xi_1+\xi_2-\xi_3-\xi_4=0} i 
b_4(\xi_1,\xi_2,\xi_3,\xi_4) \widehat{|u|^2u}(\xi_1) \hat u(\xi_2) \overline{\hat
u(\xi_3) \hat u (\xi_4)} d\sigma
\]
To achieve the cancellation of the quadrilinear form we define $b_4$ 
by 
\begin{equation}
b_4(\xi_1,\xi_2,\xi_3,\xi_4) = -  \frac{a(\xi_1)+a(\xi_2)-a(\xi_3)-a(\xi_4)}
{\xi_1^2+\xi_2^2- \xi_3^2- \xi_4^2}, \qquad (\xi_1,\xi_2,\xi_3,\xi_4)
\in P_4
\label{b4}\end{equation}
Summing up the result of our computation, we obtain
\begin{equation}
\frac{d}{dt} (E_0(u) +E_1(u)) = R_6(u)
\end{equation}
In order to estimate the size of $E_1(u)$ and of $R_6$ we need to
understand the size and regularity of $b$. A-priori $b$ is only
defined on the diagonal $P_4$. However, in order to separate variables
easier it is convenient to extend it off diagonal in a smooth way.

\begin{proposition}
  Assume that $a \in S^s_\epsilon$ with $s+\e \leq 0$. Then for each dyadic $\lambda \leq
  \alpha \leq \mu$ there is an extension of $b_4$ from the diagonal
  set
\[
 \{ (\xi_1,\xi_2,\xi_3,\xi_4) \in P_4, \ |\xi_1| \approx 
\lambda,\ |\xi_3| \approx \alpha, \ |\xi_2|,|\xi_4| \approx \mu \}
\]
to the full dyadic set 
\[
\{ \ |\xi_1| \approx 
\lambda,\ |\xi_3| \approx \alpha, \ |\xi_2|,|\xi_4| \approx \mu \}
\]
which satisfies the size and regularity conditions
\begin{equation}
| \partial_1^{\beta_1} \partial_2^{\beta_2}\partial_3^{\beta_3}
\partial_4^{\beta_4} 
 b_4(\xi_1,\xi_2,\xi_3,\xi_4)| \lesssim a(\lambda) \alpha^{-1}\mu^{-1}
 \lambda^{-\beta_1} \alpha^{-\beta_2} \mu^{-\beta_3-\beta_4}
\end{equation}
Here the implicit constants are independent of $\lambda,\alpha,\mu$
but may depend on the $\beta_j$'s.
\label{pbbounds}\end{proposition}

\begin{proof}
We first note that on $P_4$ we have the factorization 
\[
\xi_1^2+\xi_2^2 -\xi_3^2 -\xi_4^2 = 2(\xi_1 - \xi_3)(\xi_1-\xi_4)
\]
along with all versions of it due to the symmetries of $P_4$.
We consider several cases:

(a) $\lambda \ll \alpha \leq \mu$. Then the extension of $b_4$ is
defined using the formula
\[
b_4(\xi_1,\xi_2,\xi_3,\xi_4) = -  \frac{a(\xi_1)+a(\xi_2)-a(\xi_3)-a(\xi_4)}
{2(\xi_1 - \xi_3)(\xi_1-\xi_4)}
\]
and its size and regularity properties are straightforward since
$|\xi_1 - \xi_3| \approx \alpha$ and $|\xi_1 - \xi_4| \approx \mu$.

(b) $\lambda \approx \alpha \ll \mu$. Then the extension of $b_4$ is
defined using the formula
\[
b_4(\xi_1,\xi_2,\xi_3,\xi_4) = -  \frac{a(\xi_1)-a(\xi_3)}
{2(\xi_1 - \xi_3)(\xi_1-\xi_4)} -  \frac{a(\xi_2)-a(\xi_4)}
{2(\xi_4 - \xi_2)(\xi_1-\xi_4)}
\]
Now only $|\xi_1 - \xi_4| \approx \mu$ is an elliptic factor, while
the remaining quotients  exhibit suitable cancellation properties.

(c) $\lambda \approx \alpha \approx \mu$. 
Then the extension of $b_4$ is
defined by
\[
b_4(\xi_1,\xi_2,\xi_3,\xi_4) = -  \frac{a(\xi_1)+a(\xi_2)-a(\xi_1+\xi_2-\xi_4)-a(\xi_4)}
{2(\xi_1-\xi_4)(\xi_2-\xi_4)}
\]
To see that this is a smooth function on the appropriate scale we
write it in the form
\[
\begin{split}
b_4(\xi_1,\xi_2,\xi_3,\xi_4) &= \frac{1}{2(\xi_2-\xi_4)}\left(
  \frac{a(\xi_4)-a(\xi_1)}{\xi_1-\xi_4} -
  \frac{a(\xi_2)-a(\xi_2+\xi_1-\xi_4)}{\xi_1-\xi_4}\right)
\\ & = \frac{ q(\xi_4,\xi_1) - q(\xi_4+(\xi_2-\xi_4),\xi_1+(\xi_2-\xi_4))}{2(\xi_2-\xi_4)}
\end{split}
\]
where $q$ is the smooth function 
\[
q(\xi,\eta) = \frac{q(\xi) -q(\eta)}{\xi-\eta}
\]

\end{proof}

The contribution of $E_1$ to the energy is easy to control,

\begin{proposition}
Assume that $a \in S^s_\e$ with $-\frac12 < s-\e < s+\e \leq 0$.  Then 
\begin{equation}
|E_1(u)| \lesssim E_0(u)^2 
\end{equation}
\label{e1bd}\end{proposition}
We note that the threshold $s = -\frac12$ in the proposition is
consistent with the scaling.

\begin{proof}
We organize the four frequencies $\xi_1,\xi_2,\xi_3$ and $\xi_4$ 
in dyadic regions of size $\lambda \leq  \alpha \leq \mu = \mu$.  
The pointwise bound on $b$ is all we need for the proof 
since in the Fourier space one sees that only the size of the Fourier
transform matters. For a function $u$ we define $\tilde u$ by $\hat
{\tilde{u}} = |\hat u|$.  We obtain
\[
\begin{split}
  |E_1(u)| & \lesssim \sum_{\lambda \leq \alpha \leq \mu = \mu}
  |E_1 (u_\lambda, u_\alpha,u_\mu,u_\mu)|\\
  & \lesssim \sum_{\lambda \leq \alpha \leq \mu = \mu} a(\lambda)
  \alpha^{-1} \mu^{-1}
  \|\tilde u_\lambda \tilde u_\alpha \tilde u_\mu \tilde
  u_\mu\|_{L^1} \\ & \lesssim
 \sum_{\lambda \leq \alpha \leq \mu = \mu} a(\lambda)
  \alpha^{-1} \mu^{-1} \|\tilde u_\lambda \|_{L^\infty} \|\tilde u_\alpha\|_{L^\infty}
 \| \tilde u_\mu\|_{L^2}^2 
 \\ & \lesssim
 \sum_{\lambda \leq \alpha \leq \mu = \mu} \lambda^\frac12 a(\lambda)  \alpha^{-\frac12}
  \mu^{-1} \|\tilde u_\lambda\|_{L^2} \|\tilde u_\alpha\|_{L^2}
 \| \tilde u_\mu\|_{L^2}^2 
\\ & \lesssim  E_0(u)^2 \sum_{\lambda \leq \alpha \leq \mu} \frac{\lambda a(\lambda)}
{\alpha a(\alpha)  \mu^2a^2(\mu)}
\end{split}
\]
where at the last step we have used Cauchy-Schwartz with respect to
all parameters. Since $s-\e > -\frac12$ it follows that the function
$\lambda a(\lambda)$ increases polynomially with respect to $\lambda$
therefore the last sum is finite.

\end{proof}

The more difficult result we need to prove is

\begin{proposition} \label{r6} Assume that $a \in S^s_\e$ with $s+\e
  \leq 0$ and $s \geq -\frac16$. Then we have
\begin{equation}
\left|\int_0^1 R_6(u) dxdt\right| \lesssim \| u\|_{X^s}^6
\end{equation}
\end{proposition}
We note that the restriction on the symbol $a$ above is very mild.
This is because, as one can see in the proof below, there is always a
low frequency gain in the estimates. The main condition $s \geq
-\frac16$ arises in the summation with respect to high frequency
factors.

\begin{proof}
We consider a full dyadic decomposition and express the above integral
in the Fourier space as a sum of terms of the form
\[
I = \int_{0}^1 \int_{P_6} b_4(\xi_1,\xi_2,\xi_3,\xi_0) \hat
u_{\lambda_1}(\xi_1) \overline{\hat u_{\lambda_2}(\xi_2)} \hat
u_{\lambda_3}(\xi_3) P_{\lambda_0}(\overline{ \hat
  u_{\lambda_4}(\xi_4)} \hat u_{\lambda_5}(\xi_5) \overline{ \hat
  u_{\lambda_6}(\xi_6)}) d\xi dt
\]
where 
\[ 
P_6 = \{ \xi_1+\xi_3+\xi_5 = \xi_2 + \xi_4 +\xi_6 \}, \qquad \xi_0 =
\xi_1-\xi_2+\xi_3
\]
Since $b$ is smooth in each variable on the corresponding dyadic scale
we can expand it in a rapidly convergent Fourier series and separate
the variables. Hence from here on we replace $b$ by the pointwise
bound given in Proposition~\ref{pbbounds}.  There are two cases to
consider:

{\bf Case 1:}  $\lambda_0 \ll \lambda_{4,5,6}$. Then for the frequency
$\lambda_0$ factor we use Lemma~\ref{lhhh}. We denote
\[
(\lambda_1, \lambda_2, \lambda_3, \lambda_0) =
(\lambda,\alpha,\mu,\mu),
\qquad \lambda \leq \alpha \leq \mu
\]
and 
\[
f_{\lambda_0} =  \sum_{\lambda_{4,5,6} \gg \lambda_0} P_{\lambda_0}(\overline{ \hat
  u_{\lambda_4}(\xi_4)} \hat u_{\lambda_5}(\xi_5) \overline{ \hat
  u_{\lambda_6}(\xi_6)}) 
\]
We also recall the bound for $b_4$, namely 
\[
|b_4| \lesssim a(\lambda) \alpha^{-1} \mu^{-1}
\]

{\bf Case 1(a)} $\lambda_0 = \mu$. We consider the three possible 
terms in $f_{\lambda_0}$. For the $L^1 L^2$ term we bound
$u_\lambda$, $u_\alpha$ in $L^\infty$ and $u_{\mu}$ in $L^\infty_x L^2_t$.
We also sum up with respect to $\mu^{4s}$ time intervals.
This yields
\[
\begin{split}
|I| &\lesssim \mu^{-4s} \lambda^{-s} \alpha^{-s} \mu^{-s} \lambda^\frac12
\alpha^\frac12 \mu^{-1-3s} a(\lambda) \alpha^{-1} \mu^{-1} \| u\|_{X^s}^6
\\
& = \lambda^{\frac12-s} a(\lambda) \alpha^{-s-\frac12} \mu^{-2-8s}\| u\|_{X^s}^6
\end{split}
\]
The summation with respect to $\lambda$, $\mu$ and $\alpha$ is
straightforward if $s > -\frac14$.

For the $L^\frac43_x L^1_t$ term in $f$ we bound
$u_\lambda$, $u_\alpha$ in $L^\infty$ and $u_{\mu}$ in $L^4 L^\infty$.
This yields
\[
|I| \lesssim \mu^{-4s} \lambda^{-s} \alpha^{-s} \mu^{-s}
\lambda^\frac12 \alpha^\frac12 \mu^\frac14 \mu^{-\frac54-3s}
a(\lambda) \alpha^{-1} \mu^{-1} \| u\|_{X^s}^6 
\]
which gives the same outcome as in the previous case.

For  the $L^2$ part of  $f$ at modulation $\sigma \gg \mu^2$ we note 
that at least one other factor must also have modulation at least
$\sigma$. We bound that factor in $L^2$ and the other two in
$L^\infty$ to obtain
\[
\begin{split}
|I| &\lesssim \mu^{-4s} \lambda^{-s} \alpha^{-s} \mu^{-s}
 \alpha^\frac12 \mu^\frac12 \mu^{-1-3s}
a(\lambda) \alpha^{-1} \mu^{-1} \| u\|_{X^s}^6 
\\ &= a(\lambda) \lambda^{-s}
\alpha^{-s-\frac12} \mu^{-\frac32-8s}\| u\|_{X^s}^6
\end{split}
\]
which is then summed  with respect to $\lambda$, $\mu$ and $\alpha$
provided that $s > -\frac15$.

{\bf Case 1(b)} $\lambda_0 = \alpha \ll \mu$.  This case is simpler; As a
consequence of Lemma~\ref{lhhh} and of Bernstein's inequality we have
the $L^2$ type bounds
\begin{equation}
(\alpha\mu)^{-\frac12} \| Q_{\leq \alpha \mu} f_{\lambda_0}\|_{L^2} + 
 \| Q_{\gg \alpha \mu} f_{\lambda_0}\|_{X^{0,-\frac12,1}} \lesssim
 \lambda_0^{-1-3s} \|u\|_{X^s}^3 
\label{flo}\end{equation}
which is all that we need in the sequel. 

For the low modulation part of $f_{\lambda_0}$ we use the first part
of \eqref{flo}. By orthogonality we can localize the frequency $\mu$
factors to $\alpha$ intervals. Then we use the bilinear $L^2$ estimate
for $u_\lambda u_\mu$ and the pointwise bound for the other $u_\mu$
factor. This gives
\[\begin{split}
|I| &\lesssim \alpha^{-2s} \mu^{-2s} \lambda^{-s}  \mu^{-2s}
 \alpha^\frac12 \mu^{-\frac12}  \alpha^{-1-3s} (\alpha \mu)^\frac12
a(\lambda) \alpha^{-1} \mu^{-1} \| u\|_{X^s}^6 \\& = a(\lambda) \lambda^{-s}
\alpha^{-1-5s} \mu^{-1-4s}\| u\|_{X^s}^6
\end{split}
\]
The factor $\alpha^{-2s} \mu^{-2s}$ above comes from summation over
small time intervals. This is better than the earlier $\mu^{-4s}$
factor because $Q_{\leq \alpha \mu} f_{\lambda_0}$ is square
integrable on the better $\alpha^{4s}$ time scale. This is summable
with respect to $\lambda$, $\alpha$ and $\mu$ if $s \geq -\frac29$.

For the $L^2$ part of $f_{\lambda_0}$ at modulation $\sigma \gg \alpha
\mu$ we note that at least one other factor must also have modulation
at least $\sigma$. We bound that factor in $L^2$ and the other two in
$L^\infty$ to obtain
\[
\begin{split}
|I| & \lesssim \alpha^{-2s} \mu^{-2s} \lambda^{-s}  \mu^{-2s}
 \alpha  \alpha^{-1-3s}
a(\lambda)\alpha^{-1} \mu^{-1} \| u\|_{X^s}^6 
\end{split}
\]
which is the same result as above. Note that only an
$\alpha^{\frac12}$ factor is lost in the pointwise bound for $u_\mu$
due to the additional frequency localization to an interval of size
$\alpha$.

{\bf Case 1(c)} $\lambda_0 = \lambda \ll \alpha$. For the part of
$f_{\lambda_0}$ with modulation $\lesssim \alpha \mu$ we bound
$u_\alpha u_\mu$ in $L^2$ and the other $u_\mu$ in $L^\infty$. This
works even if $\alpha \approx \mu$ as two of the $\mu$ sized
frequencies must be $\mu$ separated.  We obtain
\[
\begin{split}
 |I| & \lesssim \lambda^{-2s} \mu^{-2s} \alpha^{-s}  \mu^{-2s}
 \alpha^\frac12 \mu^{-\frac12}  \lambda^{-1-3s} (\alpha \mu)^\frac12
a(\lambda) \alpha^{-1} \mu^{-1} \| u\|_{X^s}^6 \\ & = a(\lambda) \lambda^{-5s-1}
\alpha^{-s} \mu^{-1-4s}\| u\|_{X^s}^6
\end{split}
\] 
which can be summed up for $s \geq -\frac15$.

If we consider the part of $f_{\lambda_0}$ with modulation $\sigma \gg
\alpha \mu$ then another factor must have modulation at least
$\sigma$. We bound that factor in $L^2$ and the other two in
$L^\infty$ as in Case 1(b).

{\bf Case 2:} $\lambda_0 \gtrsim \min\{\lambda_4,\lambda_5,\lambda_6\}$. 
Without any restriction in generality we assume that
\[
\lambda_1 \leq \lambda_2 \leq \lambda_3, \qquad \lambda_4 \leq
\lambda_5 \leq \lambda_6
\]
Then we must have 
\[
\lambda_4 \leq \lambda_0 \leq \lambda_3
\]
We  can distribute $P_{\lambda_0}$ to each factor and also assume that 
the $\lambda_5$, $\lambda_6$ factors have frequency spread 
at most $\lambda_0$.

Denote 
\[
\{ \lambda_1, \lambda_2 ,\lambda_3, \lambda_0\} = \{ \lambda,
\alpha,\mu,\mu\}, \qquad \lambda \leq \alpha \leq \mu
\]

{\bf Case 2a:} $\lambda_0 = \mu$.

{\bf Case 2a(i):} $\lambda_5 = \lambda_6 \gg \mu$.  We use the
bilinear $L^2$ estimate for the products $u_\lambda u_{\lambda_5}$ and $u_{\lambda_4}
u_{\lambda_6}$ and the $L^\infty$ bound for $u_{\lambda}$, $u_\alpha$ and add up
with respect to $\lambda_6^{-4s}$ time intervals. We obtain
\[
\begin{split}
|I| & \ \lesssim \lambda_6^{-4s} \lambda^{-s} \alpha^{-s} \mu^{-2s}
\lambda_6^{-2s} \lambda^{\frac12} \alpha^\frac12 \lambda_6^{-1}
a(\lambda) \alpha^{-1} \mu^{-1} \prod \| u_{\lambda_j}\|_{X^s} 
\\ & \  \lesssim  
 a(\lambda) \lambda^{\frac12-s} \alpha^{-s-\frac12} 
  \mu^{-2s-1} \lambda_6^{-1-6s}\prod \| u_{\lambda_j}\|_{X^s} 
\end{split}
\]
which we sum easily with respect to the parameters $\lambda_j$ subject
to the restrictions above. We note that the summation with respect to
$\lambda_5=\lambda_6$ requires imposes the tight restriction $s \geq -\frac16$.

{\bf Case 2a(ii):} $\lambda_5 \leq \lambda_6 = \mu$ and $\alpha \ll
\mu$. Then we use the pointwise bound for $u_\lambda$, the bilinear
$L^2$ estimate for $u_\alpha u_{\mu}$ and the $L^6$ Strichartz
estimate for the remaining two factors; finally, we sum up with
respect to $\mu^{-4s}$ time intervals.  We obtain
\[
 \begin{split}
|I|
&\ \lesssim \mu^{-4s}  \lambda^{-s} \alpha^{-s} \mu^{-4s} \lambda^{\frac12}
\mu^{-\frac12} a(\lambda) \alpha^{-1} \mu^{-1} \|u_\lambda\|_{X^s}  \|u_\alpha\|_{X^s}
 \| u_\mu\|_{X^s}  \| u\|_{X^s}^3 
\\&\ \lesssim a(\lambda) \lambda^{\frac12-s}
 \alpha^{-s-1} \mu^{-\frac32-8s}  \|u_\lambda\|_{X^s}  \|u_\alpha\|_{X^s}
 \| u_\mu\|_{X^s}  \| u\|_{X^s}^3 
\end{split}
\]
The summation with respect to $\lambda$, $\alpha$ and $\mu$ requires
$s \geq \frac3{16}$.

{\bf Case 2a(iii):}  $\lambda_5 \leq \lambda_6 = \mu$ and $\alpha =
\mu$. Then we use the $L^6$ Strichartz estimate for all the factors to obtain
\[
|I| \lesssim  a(\lambda) \mu^{-4s}  \mu^{-6s} \mu^{-2} \|u\|_{X^s}^6 =
a(\lambda) \mu^{-10s-2}
\|u\|_{X^s}^6
\]

{\bf Case 2b:} $\lambda_0 = \alpha \ll \mu$.

{\bf Case 2b(i):} $\lambda_5 = \lambda_6 \gg \mu$.  Then we use the
bilinear $L^2$ estimate for $u_\mu u_{\lambda_5}$ and $u_{\mu} u_{\lambda_6}$
and $L^\infty$ for $u_{\lambda}$, $u_{\lambda_4}$ and add up with
respect to $\lambda_6^{4s}$ time intervals. We obtain exactly the same
bound as in Case 2a(i).

{\bf Case 2b(ii):} $\lambda_6 \leq \lambda_7 \leq \mu$. Then we are in
the same situation as in Case 2a(ii).

{\bf Case 2c:} $\lambda_0 = \lambda \ll \mu$. Then we can argue in the
same way as in Case 2b.

\end{proof}

The final step in the paper is to use Proposition~\ref{r6} in order to
conclude the proof of Proposition~\ref{penergy}. We have
\[
\| u_0\|_{H^s}^2 = \sum_{\lambda} \lambda^{2s} \|u_{0\lambda}\|_{L^2}^2
\]
Then the following result is straightforward:

\begin{lemma} There is a  sequence $\{ \beta_\lambda\}$ with the following
properties:

(i) $\lambda^{2s} \|u_{0\lambda}\|_{L^2}^2 \leq \beta_\lambda
\|u_0\|_{H^s}^2$.

(ii) $ \sum \beta_\lambda \lesssim 1 $.

(iii) $\beta_\lambda$ is slowly varying in the sense that
\[
|\log_2 \beta_\lambda -\log_2 \beta_\mu| \leq \frac{\epsilon}{2} 
|\log_2 \lambda- \log_2 \mu|
\]
\end{lemma}

The sequence $\beta_\lambda$ is easy to produce. One begins 
with the initial guess
\[
\beta_\lambda^0 = \frac{\lambda^{2s} \|u_{0\lambda}\|_{L^2}^2}{\|u_0\|_{H^s}^2}
\]
which satisfies (i) and (ii) but might not be slowly varying. To
achieve (iii) we mollify $\beta_\lambda^0$ on the dyadic scale and set
\[
\beta_\lambda = \sum_\mu 2^{-\frac{\e}2 |\log_2 \lambda - \log_2 \mu|}
  \beta_\mu^0 
\]

The sequence $\beta_\lambda$ will play the role of frequency localized
energy threshold.  Precisely, we assume that 
\begin{equation} \Vert u \Vert_{l^2 L^\infty L^2} \ll  1.  \label{smallenergy}
\end{equation} 
and we will show that 
\begin{equation}
\sup_t \lambda_0^{s} \|u_{\lambda_0}(t)\|_{L^2} \lesssim
\beta_{\lambda_0}^\frac12  
(\|u_0\|_{H^s} +        \|u\|_{X^s}^3)
\label{el0}\end{equation}
which by (ii) implies the conclusion of Proposition~\ref{penergy}.

In order to prove \eqref{el0} for some frequency $\lambda_0$ we define
the sequence
\[
a_{\lambda} = \lambda^{2s}  \max\{ 1, \beta_{\lambda_0}^{-1}
2^{- \epsilon |\ln \lambda- \ln \lambda_0|}\}
\]
We obtain   using the slowly varying condition (iii) 
\[
\sum_{\lambda} a(\lambda) \|u_{0\lambda}\|_{L^2}^2 \lesssim
\sum_{\lambda} \lambda^{2s} \|u_{0\lambda}\|_{L^2}^2 + 2^{-
  \frac{\epsilon}2 |\ln \lambda- \ln \lambda_0|}\lambda^{2s}
\beta_\lambda^{-1} \|u_{0\lambda}\|_{L^2}^2 \lesssim  \| u_0\|_{H^s}^2
\]

Correspondingly we find a function $a(\xi) \in S^s_\epsilon$  so that
\[
a(\xi) \approx a(\lambda), \qquad |\xi| \approx \lambda
\]
>From \eqref{smallenergy} we obtain  $ \sup_t  E_0(u(t)) \ll 1$. 
Now we use the energy estimates in Proposition~\ref{r6} for this
choice of $a$. By Proposition~\ref{e1bd} the $E_1$ component of the
energy is controlled by $E_0$, so we obtain
\[
\left(\sum_{\lambda} a(\lambda) \|u_{\lambda}(t)\|_{L^2}^2\right)^\frac12 \lesssim \|
u_0\|_{H^s} + \|u\|_{X^s}^3
\]
which at $\lambda = \lambda_0$ gives  \eqref{el0}.

\bibliographystyle{plain}
\bibliography{nls}
\end{document}